\documentclass{article}

\usepackage{proof}
\usepackage{latexsym}
\usepackage{amsfonts}
\usepackage{amssymb}
\usepackage{cite}

\newcommand{\blem}{\begin{lemma}}
\newcommand{\elem}{\end{lemma}}
\newcommand{\bth}{\begin{theorem}}
\newcommand{\ethm}{\end{theorem}}
\newcommand{\benu}{\begin{enumerate}}
\newcommand{\eenu}{\end{enumerate}}
\newcommand{\bdes}{\begin{description}}
\newcommand{\edes}{\end{description}}
\newcommand{\bdf}{\begin{definition}}
\newcommand{\edf}{\end{definition}}
\newcommand{\bcor}{\begin{cor}}
\newcommand{\ecor}{\end{cor}}
\newcommand{\bprp}{\begin{proposition}}
\newcommand{\eprp}{\end{proposition}}
\newcommand{\bmlem}{\begin{mlemma}}
\newcommand{\emlem}{\end{mlemma}}
\newcommand{\bclm}{\begin{claim}}
\newcommand{\eclm}{\end{claim}}
\newcommand{\bprf}{{\bf Proof}.\hspace{2mm}}

\newcommand{\eprf}{\hspace*{\fill} $\Box$}

\newcommand{\beqn}{\begin{equation}}
\newcommand{\eeqn}{\end{equation}}
\newcommand{\beqnarr}{\begin{eqnarray}}
\newcommand{\eeqnarr}{\end{eqnarray}}
\newcommand{\beqnarrs}{\begin{eqnarray*}}
\newcommand{\eeqnarrs}{\end{eqnarray*}}
\newcommand{\ull}{\underline{\ll}}

\newcommand{\spand}{\,\&\,}

\newtheorem{theorem}{Theorem}[section]
\newtheorem{definition}[theorem]{Definition}
\newtheorem{proposition}[theorem]{Proposition}
\newtheorem{lemma}[theorem]{Lemma}
\newtheorem{cor}[theorem]{Corollary}
\newtheorem{mlemma}[theorem]{Main Lemma}
\newtheorem{claim}[theorem]{Claim}

\newcommand{\alp}{\alpha}

\newcommand{\veps}{\varepsilon}
\newcommand{\del}{\delta}
\newcommand{\Del}{\Delta}
\newcommand{\ome}{\omega}
\newcommand{\Ome}{\Omega}
\newcommand{\bet}{\beta}
\newcommand{\gam}{\gamma}
\newcommand{\Gam}{\Gamma}
\newcommand{\kap}{\kappa}
\newcommand{\sig}{\sigma}
\newcommand{\Sig}{\Sigma}

\newcommand{\lam}{\lambda}
\newcommand{\Lam}{\Lambda}

\newcommand{\fal}{\forall}
\newcommand{\exi}{\exists}

\newcommand{\Rarw }{\Rightarrow}

\newcommand{\lrarw}{\leftrightarrow}
\newcommand{\Lrarw}{\Leftrightarrow}

\newcommand{\calh}{{\cal H}}

\newcommand{\cals}{{\cal S}}

\newcommand{\calL}{{\cal L}}

\newcommand{\msfiv}{\mbox{\hspace{5mm}}}

\newcommand{\dg}{\mbox{{\rm dg}}}

\title{Cut-elimination for {\sf SBL}}
\author{Toshiyasu Arai
\\
Graduate School of Science,
Chiba University
\\
1-33, Yayoi-cho, Inage-ku,
Chiba, 263-8522, JAPAN
\\
tosarai@faculty.chiba-u.jp
}
\date{}
\begin{document}
\maketitle
\begin{abstract}
In this paper we give
 a terminating cut-elimination procedure for a logic calculus ${\sf SBL}$.
${\sf SBL}$ corresponds to the second order arithmetic $\Pi^{1}_{2}$-Separation and Bar Induction.
\end{abstract}

\section{Introduction}\label{sec:prel3a}

Let $\Pi^{1}_{2}\mbox{-Sep}+\mbox{BI}(=\Del^{1}_{2}\mbox{-CA}+\mbox{BI})$ denote the subsystem of second order arithmetic with $\Pi^1_2$-Separation and Bar Induction. 
$\Pi^{1}_{2}\mbox{-Sep}+\mbox{BI}$ is proof-theoretically equivalent to the set theory
$\mbox{KP}i$ for recursively inaccessible universes.
K. Sch\"utte \cite{Schuette87} gives an upper bound $\psi_0I$ for the proof theoretic ordinal of 
 $\Pi^{1}_{2}\mbox{-Sep}+\mbox{BI}$. 
The ordinal $\psi_0I$ is the order type of an initial segment of the recursive notation system $T(I)$ of ordinals introduced by W. Buchholz and Sch\"utte \cite{Buchholz-Schuette83}.
G. J\"ager\cite{Jaeger83} shows the wellfoundedness up to each ordinal$<\psi_{0}I$
 in the S. Feferman's\cite{Feferman} constructive theory $T_{0}$, which is interpretable in
$\Pi^{1}_{2}\mbox{-Sep}+\mbox{BI}$.
Thus the proof-theoretic ordinal of $\Pi^{1}_{2}\mbox{-Sep}+\mbox{BI}$ and of $T_{0}$ is 
shown to be equal to $\psi_{0}I$.
J\"ager's proof is based on Ausgezeichnete Klass introduced by Buchholz\cite{Buchholz75}.

The analysis of the derivations in $\Pi^1_2\mbox{-Sep}+\mbox{BI}$ due to Sch\"utte
is based on the Buchholz's $\Ome_{\mu+1}$-rule, and the system $(T(I),<)$ is utilized indirectly:
in fact the totally defined collapsing functions $d$ and $d_{\sig}$ appear in the analysis,
which are also introduced in \cite{Buchholz-Schuette83}.

On the other side G. Takeuti \cite{Takeuti87} uses his systems of ordinal diagrams
directly for a proof theory of $\Pi^1_1$-Comprehension. 
The definition of ordinal diagrams is closely related to the cut-elimination procedure due to him. 
But unfortunately Takeuti's systems of ordinal diagrams are equipped with
many order relations and are 
too small to handle such a strong theory $\Pi^{1}_{2}\mbox{-Sep}+\mbox{BI}$.

Turning to the problem of the cut-elimination in second order, and higher order logic calculi
(known as Takeuti's Fundamental Conjecture),
W. Tait\cite{Tait66} proves the cut-eliminability (Hauptsatz)
 for the classical second order (full impredicative) logic calculus
based on the Sch\"utte's\cite{Schuette60} reformulation
of it by means of a semantical notion, \textit{semivaluation}.

Given these advances in 1980's, we had introduced a system $(O(I),<)$ of ordinal diagrams
and proved a cut-elimination theorem for a logic calculus ${\sf SBL}$
in the style of Gentzen-Takeuti \cite{Gentzen38, Takeuti87}
by transfinite induction on the system.
This was done in the original version of this paper written in 1988.
The system $(O(I),<)$ of ordinal diagrams was obtained as a kind of mixture of totally defined collapsing functions $d_{\sig}$ in \cite{Buchholz-Schuette83} and Takeuti's ordinal diagrams. 
Specifically $d_{\sig}$ is a primitive constructor of ordinal terms in $O(I)$,
whereas it is a derived term in \cite{Buchholz-Schuette83}.
In the original version of this paper it was shown that
 each initial segment determined by $\alp<\Ome_{1}\in O(I)$ is well-founded.
The proof is formalizable in $T_{0}$ as in the J\"ager's proof\cite{Jaeger83}.
 This was a starting point for us to construct larger notation systems of ordinals, e.g., in \cite{odM}.
${\sf SBL}$ corresponds to the system $\Pi^{1}_{2}\mbox{-Sep}+\mbox{BI}$ in the sense that
the Hauptsatz (normal form theorem) for ${\sf SBL}$ is equivalent to the 1-consistency of
$\Pi^{1}_{2}\mbox{-Sep}+\mbox{BI}$ over a weak theory, e.g., over $I\Sig_{1}$.
The proof of the cut-elimination in the original version was inspired from Sch\"utte's proof in \cite{Schuette87}.

In the present version let us update the original proof
 via the partially defined collapsing functions $\psi_{\sig}$
and the operator controlled derivations both due to Buchholz\cite{Buchholz92}.

In section \ref{sect:1} let us recall the collapsing functions $\psi_{\sig}$ up to $\sig\leq I$,
where $I$ is the least weakly inaccessible cardinal.
A wellfoundedness proof in $T_{0}$ is omitted in the present version
since it should not be hard.
In subsection \ref{sect:3} we define an essentially less than relation $\alp\ll\bet\,\{\eta\}$ 
for ordinals $\alp,\bet,\eta$ in terms of Skolem hulls $\calh_{\gam}(\psi_{\sig}\gam)$.
In section \ref{sect:4} a second order logic calculus ${\sf SBL}$ in introduced.
In section \ref{sect:5} we introduce a stratified logic calculus ${\sf SBL}'$ following 
Sch\"utte\cite{Schuette87}.
${\sf SBL}$ is then embedded in ${\sf SBL}'$, and a cut-free proof in ${\sf SBL}'$ 
denotes essentially a cut-free proof in ${\sf SBL}$.
For each proof $P$ in ${\sf SBL}'$ we assign an ordinal $o(P)<\psi_{\Ome_{1}}\veps_{I+1}$ in such a way that if $P$ contains a cut rule, then we can construct another proof $P'$ of the same end sequent in ${\sf SBL}'$ such that $o(P')<o(P)$ (Main Lemma \ref{mlem}).
It turns out that each proof appearing in the cut-elimination procedure
enjoys some conditions on assigned ordinals, which are 
spelled out in Definition \ref{df:31}.\ref{df:31.8}.
Restrictions similar to these conditions are found in \cite{Buchholz92}.
So our proof seems to be a finitary analogue to the proof through operator controlled derivations.

The final section \ref{sect:6} is devoted to a proof of Main Lemma \ref{mlem}.

\section{Collapsing functions $\psi_{\sig}$}\label{sect:1}

Let $I$ denote the least weakly inaccessible cardinal, and $\Ome_{\alp}:=\ome_{\alp}$ for $0<\alp<I$.
Put $\Ome_{0}:=0$ and $R=\{\Ome_{\alp+1}: \alp<I\}\cup\{I\}=\{\sig\leq I: \ome<\sig \mbox{ is regular}\}$.
$\sig,\tau,\kap$ range over elements in $R$.

In this section let us recall the collapsing functions $\psi_{\sig}\, (\sig\in R)$ due to 
W. Buchholz\cite{Buchholz92}.

\bdf
{\rm $\calh_{\alp}(X)$ denote the \textit{Skolem hull} of the set $X\cup\{0,I\}$ of ordinals under 
the functions $+, \bet\mapsto\ome^{\bet}, \bet\mapsto\Ome_{\bet}$ and 
$(\sig,\bet)\mapsto\psi_{\sig}\bet\, (\bet<\alp)$.}
\[
\psi_{\sig}\alp=
\min(\{\bet<\sig: \sig\in\calh_{\alp}(\bet) \spand \calh_{\alp}(\bet)\cap\sig\subset\bet\}\cup\{\sig\})
.\]
\edf
The following facts are shown in Lemma 4.5 of \cite{Buchholz92}.
We see that $\psi_{\sig}\alp<\sig$ from the regularity of $\sig$ and $\alp+1<\Ome_{\alp+1}$
for $\alp<I$.
When $\sig=\Ome_{\mu+1}$, we have $\mu\leq\Ome_{\mu}<\psi_{\sig}\alp<\Ome_{\mu+1}$.
Hence $\sig\in\calh_{0}(\psi_{\sig}\alp)$.
If $\alp_{0}<\alp_{1}$, then
$\calh_{\alp_{0}}(\bet)\subset\calh_{\alp_{1}}(\bet)$, and
$\psi_{\sig}\alp_{0}\leq\psi_{\sig}\alp_{1}$.
Also $\calh_{\alp}(\bet)$ is closed under the natural sum $\gam\#\del$ and the functions
$\gam\mapsto\ome^{\gam}$ and $\gam\mapsto\Ome_{\gam}$ in the reverse direction,
i.e.,
$\gam\#\del\in \calh_{\alp}(\bet) \Rarw\{\gam,\del\}\subset\calh_{\alp}(\bet)$,
$\ome^{\gam}\in\calh_{\alp}(\bet)\Rarw \gam\in\calh_{\alp}(\bet)$
and $\Ome_{\gam}\in\calh_{\alp}(\bet)\Rarw\gam\in\calh_{\alp}(\bet)$.

$\veps_{I+1}$ denotes the next epsilon number above $I$.

$\calh_{\veps_{I+1}}(0)$ is the notation system of ordinals
generated from $\{0,I\}$ by $+, \bet\mapsto\ome^{\bet}, \bet\mapsto\Ome_{\bet}$ and 
$(\sig,\bet)\mapsto\psi_{\sig}\bet\, (\bet<\veps_{I+1})$.

The computability of $\calh_{\veps_{I+1}}(0)$ together with the relation $<$ on it is seen from
the following facts.
$\psi_{\sig}\alp\in \calh_{\veps_{I+1}}(0)$ iff $\{\sig,\alp\}\subset \calh_{\veps_{I+1}}(0)\cap\calh_{\alp}(\psi_{\sig}\alp)$.

$\gam\in\calh_{\alp}(\bet) \Lrarw G_{\bet}(\gam)<\alp$,
where $G_{\bet}(0)=G_{\bet}(I)=\emptyset$,
$G_{\bet}(\alp_{0}+\cdots+\alp_{n})=\bigcup\{G_{\bet}(\alp_{i}): i\leq n\}$,
$G_{\bet}(\ome^{\alp})=G_{\bet}(\Ome_{\alp})=G_{\bet}(\alp)$,
\[
G_{\bet}(\psi_{\sig}\alp)=\left\{
\begin{array}{ll}
\emptyset & \mbox{if } \psi_{\sig}\alp<\bet
\\
G_{\bet}(\sig)\cup G_{\bet}(\alp)\cup\{\alp\}
\end{array}
\right.
\]
\benu
\item
$Fx:=\{\alp\in \calh_{\veps_{I+1}}(0): \Ome_{\alp}=\alp>0\}=\{\psi_{I}\alp: \alp\in \calh_{\veps_{I+1}}(0), \alp\in\calh_{\alp}(\psi_{I}\alp)\}\cup\{I\}$.
\item
$\Ome_{\alp}<\psi_{\Ome_{\alp+1}}\bet<\Ome_{\alp+1}$.
\item
$\psi_{\Ome_{\alp+1}}\bet<\psi_{I}\gam\Lrarw \alp<\psi_{I}\gam$.
\item
$\psi_{I}\alp<\psi_{I}\bet \Lrarw \alp<\bet$.
\eenu

In what follows $\alp,\bet,\gam,\del,\ldots$ range over ordinal terms in $\calh_{\veps_{I+1}}(0)$, and
 $\sig,\tau,\ldots$ over elements in the set 
$R=\{I\}\cup\{\Ome_{\mu+1}: \mu\in \calh_{\veps_{I+1}}(0)\}$.

\subsection{Relations $\alp\ll\bet\,\{\eta\}$}\label{sect:3}
In this subsection an essentially less than relation $\alp\ll\bet\,\{\eta\}$ is defined through
Skolem hulls $\calh_{\gam}(\psi_{\sig}\gam)$.

\bdf\label{df:ll}
{\rm For ordinal terms $\del_{0},\del_{1},\eta\in \calh_{\veps_{I+1}}(0)$,
\beqnarrs
\del_{0} \ll \del_{1}\, \{\eta\} & :\Lrarw & \del_{0}<\del_{1} \land
\fal\sig\fal\alp[\{\del_{1},\eta\}\subset\calh_{\alp}(\psi_{\sig}\alp)\Rarw 
\del_{0}\in\calh_{\alp}(\psi_{\sig}\alp)]
\\
\del_{0} \ull \del_{1}\, \{\eta\} & :\Lrarw & \del_{0}=\del_{1} \lor (\del_{0} \ll \del_{1}\, \{\eta\})
\\
\del_{0} \ll \del_{1} & :\Lrarw & \del_{0} \ll \del_{1}\, \{0\}
\\
\del_{0} \ull \del_{1} & :\Lrarw & \del_{0}=\del_{1} \lor \del_{0} \ll \del_{1} 
\eeqnarrs
}
\edf

\bprp\label{prp:ll}
\benu
\item\label{prp:ll.1}
$\del_{0}\ll\del_{1}\, \{\eta\}  \Rarw \del_{0}\#\alp\ll\del_{1}\#\alp\, \{\eta\} $.
\item\label{prp:ll.2}
$\del_{0}\ll\del_{1}\, \{\eta\}  \Rarw \ome^{\del_{0}}\ll\ome^{\del_{1}}\, \{\eta\}$.
\item\label{prp:ll.3}
Assume $\{\gam,\del_{1},\eta\}\subset\calh_{\gam}(\psi_{\sig}(\gam\#\del_{1}))$.
Then 
$\del_{0}\ll\del_{1}\, \{\eta\}    \Rarw
 \psi_{\sig}(\gam\#\del_{0})\ll\psi_{\sig}(\gam\#\del_{1})\, \{\eta\}$.
\item\label{prp:ll.4}
Assume $\gam\in\calh_{\gam}(\psi_{\Ome_{\mu+1}}\gam)$. Then
$\alp\ll\gam \spand \alp\leq\Ome_{\mu} \Rarw \alp\ll\psi_{\Ome_{\mu+1}}\gam$.
\item\label{prp:ll.5}
Assume $\alp\ll\psi_{\tau}\gam\,\{\eta\}$, $\gam\in\calh_{\gam}(\psi_{\tau}\gam)$ and 
$\fal\bet>\gam\fal\sig[\{\gam,\tau\}\subset\calh_{\bet}(\psi_{\sig}\bet) \Rarw \eta\in\calh_{\bet}(\psi_{\sig}\bet)]$.
Then $\alp\ll\psi_{\tau}\gam$.
\item\label{prp:ll.6}
Assume $\del_{0}\ll\del_{1}\,\{\eta\}$, 
$\psi_{\tau}(\gam\#\del_{0})<\psi_{\tau}(\gam\#\del_{1})$,
and $\{\psi_{\tau}(\gam\#\del_{1}),\eta\}\subset\calh_{\alp}(\psi_{\sig}\alp)$.
Then $\psi_{\tau}(\gam\#\del_{0})\in\calh_{\alp}(\psi_{\sig}\alp)$.
\eenu
\eprp
\bprf
\ref{prp:ll}.\ref{prp:ll.3}.
Assume $\{\gam,\del_{1},\eta\}\subset\calh_{\gam}(\psi_{\sig}(\gam\#\del_{1}))$ and 
$\del_{0}\ll\del_{1}\, \{\eta\}$.
Then $\del_{0}\in\calh_{\gam}(\psi_{\sig}(\gam\#\del_{1}))$, and 
$\{\gam,\del_{0},\sig\}\subset\calh_{\gam\#\del_{1}}(\psi_{\sig}(\gam\#\del_{1}))$.
Hence 
$\psi_{\sig}(\gam\#\del_{0})\in\calh_{\gam\#\del_{1}}(\psi_{\sig}(\gam\#\del_{1}))\cap\sig=\psi_{\sig}(\gam\#\del_{1})$.

Next suppose $\{\psi_{\sig}(\gam\#\del_{1}),\eta\}\subset\calh_{\alp}(\psi_{\tau}\alp)$.
We show $\psi_{\sig}(\gam\#\del_{0})\in\calh_{\alp}(\psi_{\tau}\alp)$.
We can assume $\psi_{\sig}(\gam\#\del_{1})>\psi_{\tau}\alp$.
Then $\{\sig,\gam,\del_{1},\eta\}\subset\calh_{\alp}(\psi_{\tau}\alp)$ and $\gam\#\del_{1}<\alp$.
Therefore $\{\sig,\gam,\del_{0}\}\subset\calh_{\alp}(\psi_{\tau}\alp)$ and $\gam\#\del_{0}<\alp$.
These yield $\psi_{\sig}(\gam\#\del_{0})\in\calh_{\alp}(\psi_{\tau}\alp)$.
\\
\ref{prp:ll}.\ref{prp:ll.6} is seen as in Proposition \ref{prp:ll}.\ref{prp:ll.3}.
\eprf
\\

 Let $Var=\{U,V,\ldots\}$ be a countable set of (unary) second-order free variables, and
 $Var':=\{U^{\eta}: U\in Var, \eta\in Fx\}$.
Also $\Sig=\{0,I,+,\ome,\Ome,\psi\}\cup Var'$.
$Var(t)$ denotes the set of variables occurring in $t\in\Sig^{*}$(the set of finite sequences over symbols $\Sig$).

\bdf\label{df:17}
{\rm Let $\max$ be a symbol not in $\Sig$.}
\benu
\item
{\rm $\cals\subset(\Sig\cup\{\max\})^{*}$ and ordinals $od(s)\leq I$ for $s\in\cals$ are defined recursively.}
 \benu
 \item
 {\rm Each $U^{\eta}\in Var'$ is in $\cals$.
 $od(U^{\eta})=\eta$.}
 \item
 $Fx\cup\{0\}\subset\cals$.
 $od(\eta)=\eta$ {\rm for} $\eta\in (Fx\cap I)\cup\{0\}$.
 \item
 $s\in\cals \Rarw s+1\in\cals$.
 $od(s+1)=\min(od(s)+1,I)$.
 \item
 $s_{1},s_{2}\in\cals \Rarw \max(s_{1},s_{2})\in\cals$.
 \\
 $od(\max(s_{1},s_{2}))=\min(I,\max(od(s_{1}),od(s_{2})))$.
 \eenu

\item
{\rm For $s\in\cals$, a finite non-empty set $I(s)\subset Fx\cup\{0\}$ is defined recursively.}
 \benu
 \item
 $I(U^{\eta})=\{\eta\}$.
 \item
 $I(s)=\{s\}$ {\rm for} $s\in Fx\cup\{0\}$.
 \item
 $I(s+ 1)=I(s)$.
 \item
 $I(\max(s_{1},s_{2}))=I(s_{1})\cup I(s_{2})$.
 \eenu



\eenu
\edf

Note that $od(s)=I$ iff  a free variable $U^{I}$ with index $I$ occurs in $s$.
In particular if no free variable occurs in $s$, then $od(s)<I$.

\section{The logic calculus ${\sf SBL}$}\label{sect:4}

In this section a second-order logic calculus {\sf SBL} is introduced.
$\mathcal{L}$ denotes a second-order language consisting of
logical symbols $\lor,\land,\exi,\fal$,
individual constants $c,\ldots$, function symbols $f,\ldots$,
 first-order free variables $a,b,\ldots$,
first-order bound variables $x,y,\ldots$,
relation symbols $R,\ldots$,
second-order free variables $U,V,\ldots$, and
second-order bound variables $X,Y,\ldots$.
Let us assume that each relation symbol and each second-order (free or bound) variable
is \textit{unary} for simplicity, and
$\mathcal{L}$ contains an individual constant $c$ and a (unary) relation symbol $R$,
cf.\,pure variable condition in Definition \ref{df:27}.

$T$ stands for either a (unary) second-order free variable or a predicate constant.
$Tt$ and $\lnot Tt$ are \textit{prime formulas} (literals) for terms $t$.
Formulas are generated from literals by means of $\lor,\land$ and first-order and second-order
quantifications $\exi,\fal$ as usual.
Negations $\lnot A$ of formulas $A$ are defined recursively through de Morgan's law and the elimination of
double negations $\lnot(\lnot Tt):\equiv (Tt)$.

For formal expressions $E,s,t$ such as terms and proofs 
$E[s/t]$ denotes the expression obtained from $E$ by replacing some occurrences of
the expression $t$ in $E$ by the expression $s$.
Let $F$ be a formula\footnote{Strictly speaking, we should say that $F$ is a semi-formula as in \cite{Takeuti87}, but for simplicity let us call semi-formulas and semi-terms with bound variables
as formulas and terms, resp.}
 with a second-order bound variable $X$,
and $A$ a formula with a variable $x$.
Then $F[A/X]$ denotes the formula obtained from $F$ by replacing each occurrence
of $Xt$ by $A[t/x]$ and each occurrence of $\lnot Xt$ by $\lnot A[t/x]$.

\bdf\label{df:18}
\benu

\item
{\rm For formulas $A$, $VT(A)$ denotes the set of second-order free variables occurring
in a scope of a second-order quantifier in $A$.
$U\in VT(A)$ iff $U$ is tied by a second-order quantifier in $A$ in the sense of \cite{Takeuti87}.}

\item
{\rm Let $A$ be a formula.}
\beqnarrs
\fal X\, A[X/U]\in\Pi^{1}_{2} & \Lrarw & A\in\Pi^{1}_{2}
\\
\exi X\, A[X/U]\in\Sig^{1}_{2} & \Lrarw & A\in\Sig^{1}_{2}
\\
\fal X\, A[X/U]\in\Sig^{1}_{2} & \Lrarw & A\in\Pi^{1}_{2}\cap\Sig^{1}_{2} \spand U\not\in VT(A)
\\
\exi X\, A[X/U]\in\Pi^{1}_{2} & \Lrarw & A\in\Pi^{1}_{2}\cap\Sig^{1}_{2} \spand U\not\in VT(A)
\eeqnarrs
$A\in\Pi^{1}_{2}\cap\Sig^{1}_{2}$ {\rm iff $A$ is isolated in the sense of \cite{Takeuti87}.}

\item
{\rm 
An occurrence of a second-order quantifier $QX$ in a formula $QX A[X/U]$ is said to be
\textit{distinguished} if either $Q=\fal$, $\fal X A[X/U]\in\Pi^{1}_{2}$ and $U\not\in VT(A)$,
or $Q=\exi$, $\exi X A[X/U]\in\Sig^{1}_{2}$ and $U\not\in VT(A)$.
}
\eenu
\edf

\bdf\label{df:19}
{\rm The logic calculus ${\sf SBL}$.}

Axioms {\rm or} initial sequents {\rm are}
\[
(Ax) \msfiv \Gam,\lnot L,L \msfiv\mbox{{\rm for prime }} L
\]
Inference rules {\rm are the followings.}
\[
\infer[(\lor)]{\Gam}{\Gam,A_{i}}
\msfiv
\infer[(\land)]{\Gam}
{
 \Gam,A_{0}
 &
 \Gam,A_{1}
 }
\]
{\rm where in the rule $(\lor)$,
$A_{0}\lor A_{1}$ is the \textit{main formula} of $(\lor)$ and is in the \textit{lower sequent}
$\Gam$.
The formula $A_{i}\,(i=0,1)$ in the \textit{upper sequent} is the \textit{minor formula} of the $(\lor)$.
In the rule $(\land)$, $A_{0}\land A_{1}$ is the \textit{main formula},
 and is in the lower sequent $\Gam$.
Formulas $A_{i},(i=0,1)$ in the upper sequents are the minor formulas of the $(\land)$.}

\[
\infer[(\exi_{1})]{\Gam}{\Gam,F[t/x]}
\msfiv
\infer[(\fal_{1})]{\Gam}
{
 \Gam, F[a/x]
 }
\]
{\rm where in the rule $(\exi)_{1}$,
$\exi x\, F$ is the main formula of $(\exi)_{1}$ and is in the lower sequent
$\Gam$.
The formula $F[t/x]$ in the upper sequent is the minor formula of the $(\exi_{1})$.
In the rule $(\fal_{1})$, $\fal x\, F$ is the main formula,
 and is in the lower sequent $\Gam$.
The $F[a/x]$ in the upper sequent is the minor formula of the $(\fal_{1})$,
and the free variable $a$ is the \textit{eigenvariable} of the $(\fal_{1)}$,
which does not occur in the lower sequent $\Gam$.}

\[
\infer[(cut)]{\Gam,\Del}
{\Gam,\lnot C
&
C,\Del
}
\]
{\rm $C$ is the \textit{cut formula} of the $(cut)$.}

\[
\infer[(\exi_{2})]{\Gam}{\Gam,F[T/X]}
\msfiv
\infer[(\fal_{2})]{\Gam}{\Gam,F[U/X]}
\]
{\rm where in the rule $(\exi)_{2}$,
$\exi X\, F$ is the main formula of $(\exi)_{2}$ and is in the lower sequent
$\Gam$.
The formula $F[T/X]$ in the upper sequent is the minor formula of the $(\exi_{2})$.
In the rule $(\fal_{2})$, $\fal X\, F$ is the main formula,
 and is in the lower sequent $\Gam$.
The $F[U/X]$ in the upper sequent is the minor formula of the $(\fal_{2})$,
and the free variable $U$ is the \textit{eigenvariable} of the $(\fal_{2})$,
which does not occur in the lower sequent $\Gam$.}

\[
\infer[(BI)]{\Gam}{\Gam,F[A/X]}
\]
{\rm where $\exi X\, F$ is the main formula of the $(BI)$ and is in the lower sequent
$\Gam$, $F[A/X]$ is the minor formula, and
either 
\benu
\item[$(BI)_{1}$]
$\exi X\, F\in\Pi^{1}_{2}\cap\Sig^{1}_{2}$, or 
\item[$(BI)_{2}$]
$A\in\Pi^{1}_{2}\cap\Sig^{1}_{2}$.
\eenu
\[
\infer[(\Pi^{1}_{2}\mbox{{\rm -Sep}})]{\Gam}{\Gam,A\subset B}
\]
where $\exi X(A\subset X\subset B)$  is the main formula the $(\Pi^{1}_{2}\mbox{{\rm -Sep}})$ and is in the lower sequent
$\Gam$, $A\subset B$ is the minor formula, and
$A\in\Pi^{1}_{2}$ and $B\in\Sig^{1}_{2}$.}
\edf

\bth\label{th:main}{\rm Cut-elimination theorem for ${\sf SBL}$.}

There is a rewriting procedure $r$ on derivations in {\sf SBL} such that
for any {\sf SBL}-derivation $P$ of a sequent, if $P$ contains a $(cut)$,
then $r(P)$ is an {\sf SBL}-derivation of the same sequent,
and there is an $n$ such that its $n$-th iterate $r^{(n)}(P)$ is cut-free.
\end{theorem}

\bdf\label{df:20}
{\rm
\benu
\item
A formula is said to be \textit{first-order} if no second-order quantifier occurs in it.
\item
A sequent is \textit{first-order} if every formula in it is first-order.
\item
We say that the \textit{cut-elimination theorem holds for derivations in ${\sf SBL}$ ending with first-order sequents} if there is a rewriting procedure $r$ on derivations of first-order sequents
for which Theorem \ref{th:main} holds.
\eenu
}
\edf

\bprp\label{prp:8}
If the cut-elimination theorem holds for derivations in ${\sf SBL}$ ending with first-order sequents,
then the cut-elimination theorem \ref{th:main} holds for ${\sf SBL}$.
\eprp
\bprf
This is seen by cut-elimination by absorption combined with the joker translation due to
P. P\"appinghaus\cite{Paeppinghaus}.
Note that 
\[
\infer{\Gam, \lnot B_{0}}
{
\Gam,\lnot C
&
\Gam,C
}
\msfiv
\infer{\Gam,\exi Y F(Y)}
{
\Gam, F(\{x: B_{0}\})
}
\]
are admissible rules in the presence of the inference rules $(BI)_{1}$ and $(BI)_{2}$
for $B_{0}\equiv(\fal X\fal y[X(y)\to X(y)])\in\Pi^{1}_{2}\cap\Sig^{1}_{2}$, 
cf.\,Theorem 1.3 and Lemma 1.5(ii) in \cite{Paeppinghaus}, resp.
\eprf
\\

\noindent
{\bf Remark}.
Let ${\sf SBL}_{1}$ denote temporarily the calculus ${\sf SBL}$ without the rule $(BI)_{2}$.
Namely in ${\sf SBL}_{1}$, the rule $(BI)$ is restricted to the case when
the main formula $\exi X F\in\Pi^{1}_{2}\cap\Sig^{1}_{2}$.
${\sf SBL}_{1}$ is equivalent to ${\sf SBL}$ with respect to derivable sequents since
$\exi X\fal y[X(y)\lrarw A(y)]\in\Pi^{1}_{2}\cap\Sig^{1}_{2}$ for $A\in\Pi^{1}_{2}\cap\Sig^{1}_{2}$.
The reason why we introduce the superfluous $(BI)_{2}$ in ${\sf SBL}$ is as follows:
when we operate our cut-elimination procedure to an ${\sf SBL}_{1}$-derivation, then we obtain a cut-free ${\sf SBL}$-derivation with rules $(BI)_{2}$
since we need `infer $\exi Y F(Y)$ from 
$F(\{x: B_{0}\})$
in replacing the joker $J_{0}$ by $B_{0}$
in Lemma 1.5(ii) of \cite{Paeppinghaus}.
In other words, we don't have an `inner' cut-elimination theorem for ${\sf SBL}_{1}$.
It is open for us whether or not the `inner' cut-elimination theorem for ${\sf SBL}_{1}$
holds besides its intrinsic interests.


\section{The stratified logic calculus ${\sf SBL}'$}\label{sect:5}

A stratified (in German: geschichtet) calculus ${\sf SBL}'$ is introduced.

\bdf\label{df:21}
{\rm A stratified language $\calL'$ is obtained from a second-order language $\calL$
by modifying relation symbols and second-order variables as follows.}
\benu
\item
{\rm (unary) relation symbols with \textit{index} $0$:
$R^{0}$.}
\item
{\rm (unary second-order) \textit{unstratified bound} variables: $X,Y,\ldots$.}
\item
stratified {\rm variables :}
 \benu
 \item
 {\rm free variables with \textit{index} $s$ : $U^{s}$ for $s\in \cals$ and free variables $U$ in $\calL$.}
 \item
 {\rm bound variables with \textit{index} $\eta$ : $X^{\eta}$ for $\eta\in Fx$
 and bound variables $X$ in $\calL$.}
 \eenu
\eenu
\edf

When $T$ denotes a predicate constant $R$, $T^{s}:\equiv R^{0}$, i.e., $s=0$.

\bdf\label{df:22}
$\calL'$-formula $A$ {\rm is obtained from an $\calL$-formula $A$ by attaching indices as follows.}
\benu
\item
{\rm Attach the index $0$ to each predicate constant $R$ occurring in $A$.}

\item
{\rm Attach an index $s\in\cals$ to every occurrence of each free variable $U$
occurring in $A$.
The indices may depend on free variables.}

\item
{\rm Attach an index $\eta\in Fx$ to \textit{all undistinguished} quantifiers.
In a formula each undistinguished quantifier receives the same index.
Also leave distinguished quantifiers without indices.}
\eenu
\edf
$A^{\calL}$ denotes the (unstratified) $\calL$-formula obtained from an $\calL'$-formula
$A$ by erasing all indices.
Conversely $A'$ denotes ambiguously an $\calL'$-formula obtained from an $\calL$-formula $A$
by attaching some indices.

$A'\in\Pi^{1}_{2}[\Sig^{1}_{2}]$ iff $A\in\Pi^{1}_{2}[\Sig^{1}_{2}]$, resp.
$\fal^{\eta}, \exi^{\eta}$ denote stratified quantifiers $\fal X^{\eta},\exi X^{\eta}$ for a bound variable $X$.

\bdf\label{df:23}
{\rm $Gr(A),gr(A)<\ome$ for $\calL'$-formulas $A$.}
\benu
\item
 \benu
 \item
 $Gr(A)=0$ {\rm if neither $\fal^{I}$ nor $\exi^{I}$ occurs in $A$.
 
 In what follows assume that either $\fal^{I}$ or $\exi^{I}$ occurs in $A$.}
 
 \item
 $Gr(A)=\max\{Gr(A_{0}),Gr(A_{1})\}+1$ {\rm if} $A\in\{A_{0}\lor A_{1},A_{0}\land A_{1}\}$.
 \item
 $Gr(A)=Gr(B)+1$ {\rm if} $A\in\{\fal x\, B[x/u],\exi x\, B[x/u]\}$.
 \item
 $Gr(A)=1$ {\rm if} $A\in\{\fal X F,\exi X F\}$.
 \item
 $Gr(A)=\max\{2,Gr(F[R^{0}/X])+1\}$ {\rm if} $A\in\{\fal X^{I}F,\exi X^{I}F\}$.
 \eenu
\item
 \benu
 \item
 $gr(A)=0$ {\rm if $A$ is either a prime formula or of the form $QX F$.}
 \item
 $gr(A_{0}\lor A_{1})=gr(A_{0}\land A_{1})=\max\{gr(A_{0}),gr(A_{1})\}+1$.
 \item
 $gr(\exi x\, B[x/u])=gr(\fal x\, B[x/u])=gr(B)+1$.
 \item
 $gr(QX^{\eta}F)=gr(F[R^{0}/X])+1$.
 \eenu
\eenu

\edf

\bdf\label{df:24}
{\rm Let $A$ be an $\calL'$-formula.}
\benu
\item
{\rm $A\in\Sig^{I}$ (in Sch\"utte's terminology `$A$ ist klein') if $\fal^{I}$ does not occur in $A$.}
\item
{\rm An occurrence of a free variable $U^{\eta}\in Var'$ in $A$ is said to be}
 \benu
 \item
 in an index {\rm if the occurrence is in an index of a stratified (free) variable, and}
 \item
 an occurrence as a part of a formula {\rm otherwise.}
 \eenu
 
\item
{\rm $A$ is said to be \textit{stratified} if for each index $s\in\cals$ of a free variable $U^{s}$ 
occurring as a part of the formula $A$, $Var(s)=\emptyset$ and $od(s)<I$.}
\eenu
\edf

\bdf\label{df:25}
{\rm For $\calL'$-formulas $A$,
$st_{\Pi}(A)\in\cals$ if $A\in\Pi^{1}_{2}$, and $st_{\Sig}(A)\in\cals$ if $A\in\Sig^{1}_{2}$ are defined. Let $\Lam\in\{\Pi,\Sig\}$.}
\benu
\item
$st_{\Lam}(T^{s}t)=st_{\Lam}(\lnot T^{s}t)=s$.
\item
$st_{\Lam}(A_{0}\circ A_{1})=\max(st_{\Lam}(A_{0}), st_{\Lam}(A_{1}))$ {\rm for} $\circ\in\{\land,\lor\}$.
\item
$st_{\Lam}(Qx\, B[x/u])=st_{\Lam}(B)$ {\rm for} $Q\in\{\fal,\exi\}$.

{\rm In what follows let $F_{0}\equiv F[R^{0}/X]$.}
\item
\beqnarrs
st_{\Pi}(\fal X F)=st_{\Pi}(F_{0}) & ; & st_{\Sig}(\fal X F)=st_{\Pi}(F_{0})+1
\\
st_{\Sig}(\exi X F)=st_{\Sig}(F_{0}) & ; & st_{\Pi}(\exi X F)=st_{\Sig}(F_{0})+1
\eeqnarrs

\item
\beqnarrs
st_{\Pi}(\fal X^{\eta}F) & = & \max(\eta, st_{\Pi}(F_{0}))
\\
st_{\Sig}(\exi X^{\eta}F) & = & \max(\eta, st_{\Sig}(F_{0}))
\eeqnarrs
\eenu

{\rm Let $A$ be an $\calL'$-formula.}
{\rm For a variable $U^{\eta}\in Var'$ and $s\in\cals$ let $A^{[s/U^{\eta}]}$ denote the $\calL'$-formula obtained from $A$ by replacing every occurrence of $U^{\eta}$ in an index by $s$. 
$\Del^{[s/U^{\eta}]}=\{A^{[s/U^{\eta}]}:A\in\Del\}$ for sequents $\Del$, and
$P^{[s/U^{\eta}]}$ the tree of sequents obtained from a preproof $P$ by 
replacing each sequent $\Del$
in $P$ by $\Del^{[s/U^{\eta}]}$.}
\edf

\bprp\label{prp:9}
{\rm Let $A$ be an $\calL'$-formula such that $A\in\Pi^{1}_{2}$.}
 \benu
 \item\label{prp:9b1}
 $st_{\Sig}(\lnot A)=st_{\Pi}(A)$.
 \item\label{prp:9b2}
 $Var(st_{\Pi}(A))=Var(A)$.
 
{\rm In what follows assume that $A$ is stratified.}
 \item\label{prp:9b3}
 {\rm Let $\eta$ denotes the index of an undistinguished quantifier in $A$ if such a quantifier occurs.
 Otherwise let $\eta=0$.
 Also let $\nu=\max I(A)$.
 Then there is a $k<\ome$ such that $st_{\Pi}(A)=\max\{\eta,\nu+k\}$.}
 
 \item\label{prp:9b4}
 {\rm Let $A\equiv(\fal X F)$, and $U$ a variable not occurring in $A$.
 Then $st_{\Pi}(F[U^{s}/X])=st_{\Pi}(\fal X F)$ for $s=st_{\Pi}(\fal X F)$.}
 
 \eenu
 
\eprp

\bdf\label{df:26}{\rm Axioms and inference rules in ${\sf SBL}'$.}

\[
(Ax) \msfiv \Gam,\lnot A, A \mbox{   {\rm with }} Gr(A)=0
\]

{\rm $(\land),(\lor),(\fal_{1}),(\exi_{1}),(cut)$ are the same as in ${\sf SBL}$.}

\[
\infer[(th)]{\Gam,\Del}{\Gam}
\]
{\rm is the} \textit{thinning}.

\benu
\item\label{df:26.c}
critical rule.
\[
\infer[(c)]{\Gam,\exi X^{\eta}F}{\Gam,\exi X^{\eta}F,F(T^{s})}
\]
{\rm where $\eta\neq I \Rarw Var(s)=\emptyset \spand od(s)<\eta$.
$s$ the \textit{index}, and $\eta$ \textit{type} of the inference.}

\item\label{df:26.2}
distinguished rules.
 \benu
 \item\label{df:26.d1}
\[
\infer[(d1)]{\Gam,\exi X(A\subset X\subset B)}{\Gam,\exi X(A\subset X\subset B),A\subset B}
\]
{\rm where $Gr(\exi X(A\subset X\subset B))\neq 0$.}


 \item\label{df:26.d2}
\[
\infer[(d2)]{\Gam,\exi X F}{\Gam, \exi X F,F(T^{s})}
\]
{\rm where} $Gr(\exi X F)\neq 0$.
 
 \eenu
 
\item\label{df:26.BI}
\[
\infer[(BI)]{\Gam,\exi X F}{\Gam,\exi X F,F(A)}
\]
{\rm where}
 \benu
 \item\label{df:26.BI.1}
 {\rm if $\exi X F\not\in\Pi^{1}_{2}$, then $\exi X F$ is stratified, and}
 \item\label{df:26.BI.2}
 $Gr(\exi X F)=0$.
 \eenu
 
\item\label{df:26.s}
strong rules.
 \benu
 \item\label{df:26.s1}
 \[
 \infer[(s1)]{\Gam,\fal X F}{\Gam,\fal X F,F(U^{s})}
 \]
 {\rm where $Gr(\fal X F)\neq 0$, $U$ does not occur in the lower sequent, and
 $s$ is obtained from $st_{\Pi}(\fal X F)$ 
 by replacing occurrences of $I$ corresponding to an undistinguished quantifier in $\fal X F$ by $U^{I}$.
$st_{\Pi}(\fal X F)=s[I/U^{I}]$.
 $(s1)$ is of \textit{type} $I$.}
 
 \item\label{df:26.s2}
 \[
 \infer[(s2)]{\Gam,\fal X^{\eta}F}{\Gam,\fal X^{\eta} F,F(U^{U^{\eta}})}
 \]
  {\rm where $U$ does not occur in the lower sequent. $(s2)$ is of \textit{type} $\eta$.}
 \eenu

\item\label{df:26.w}
weak rule.
\[
\infer[(w)]{\Gam,\fal X F}{\Gam, \fal X F,F(U^{s})}
\]
{\rm where}
 \benu
 \item\label{df:26.w1}
 {\rm if $\fal X F\not\in\Sig^{1}_{2}$, then $\fal X F$ is stratified.}
 \item\label{df:26.w2}
 $Gr(\fal X F)=0$.
 \item\label{df:26.w3}
 $s=st_{\Pi}(\fal X F)$.
 \item\label{df:26.w4}
 {\rm $U$ does not occur in the lower sequent.}
 \eenu
$s$ {\rm is the \textit{index} of the $(w)$.}

\item\label{df:26.sub}
substitution of level $s$.
\[
\infer[(sub)^{s}]{\Gam[A/U^{s}]}{\Gam}
\]
{\rm where}
 \benu
 \item\label{df:26.sub1}
 $Var(s)=\emptyset$ {\rm and} $od(s)<I$.
 \item\label{df:26.sub2}
 {\rm $A$ is a stratified formula such that $B[A/U^{s}]$ is an $\calL'$-formula for $B\in\Gam$.}
 \item\label{df:26.sub3}
 {\rm $U$ does not occur in the lower sequent.}
 \item\label{df:26.sub4}
 {\rm any $B\in\Gam$ enjoys the followings.}
 \benu
 \item\label{df:26.sub4.1}
 {\rm $B$ is a stratified $\Pi^{1}_{2}$-formula such that $st_{\Pi}(B)\leq od(s)$.}
 \item\label{df:26.sub4.2}
 $U\not\in VT(B^{\calL})$.
 \eenu
 
 \eenu

\item\label{df:26.falred}
$\fal^{I}$-reduction of type $\eta<I$.
\[
\infer[(\fal^{I}\mbox{-red})^{\eta}]{\Del_{0}[\fal^{\eta}/\fal^{I}],\Del_{1}}{\Del_{0},\Del_{1}}
\]
{\rm where each $A\in\Del_{0}\cup\Del_{1}$ is either $A\in\Sig^{I}$ or $\lnot A\in\Sig^{I}$,
and $A[\fal^{\eta}/\fal^{I}]$ denotes the $\calL'$-formula obtained from $A$ by replacing $\fal X^{I}$ by $\fal X^{\eta}$.}

\item\label{df:26.exired}
$\exi^{I}$-reduction of type $\eta<I$.
\[
\infer[(\exi^{I}\mbox{-red})^{\eta}]{\Del_{0}[\exi^{\eta}/\exi^{I}],\Del_{1}}{\Del_{0},\Del_{1}}
\]
{\rm where $\Del_{0}\cup\Del_{1}\subset\Sig^{I}$ and
$A[\exi^{\eta}/\exi^{I}]$ denotes the $\calL'$-formula obtained from $A$ by replacing $\exi X^{I}$ by $\exi X^{\eta}$.}

\eenu

\edf

Inference rules without main formulas are $(cut),(th), (sub), (\fal^{I}\mbox{-red})$ and 
$(\exi^{I}\mbox{-red})$.

\bdf\label{df:27}
{\rm
A \textit{preproof} is a finite tree with $(Ax)$ and inference rules in ${\sf SBL}'$.

A preproof $P$ enjoys the \textit{pure variable condition} if all eigenvariables are distinct each other, each eigenvariable does not occur in the end-sequent of $P$ and
if a free variable occurs in an upper sequent of a rule, but not in the lower sequent,
then the variable is the eigenvariable of the rule.
}
\edf

Let $P$ be a preproof with the pure variable condition, and $U^{s}$ a second-order free variable occurring in $P$.
Then either the stratified variable $U^{s}$ occurs in the end-sequent, or an eigenvariable
of one of rules $(s),(w),(sub)\, J$.
Consider the latter case, and 
let $V^{\eta}$ be a variable with index $\eta$ occurring in the index $s$ of $U^{s}$.
Then the rule $J$ is either an $(s)$ or a $(w)$.
When $J$ is either an $(s1)$ or a $(w)$ with its main formula $\fal X F$, 
then either $s[I/U^{I}]=st_{\Pi}(\fal X F)$ or
$s=st_{\Pi}(\fal X F)$, cf.\,Definitions \ref{df:26}.\ref{df:26.s1} and \ref{df:26}.\ref{df:26.w}.
Hence either $V^{\eta}\equiv U^{I}$ corresponds to an undistinguished quantifier in $\fal X F$, or
$V^{\eta}$ occurs in the index of a variable occurring in the main formula $\fal X F$.
Arguing inductively,
this means that either the variable $V^{\eta}$ occurs in the end-sequent, or 
corresponds to an undistinguished quantifier in a main formula of an $(s1)$, or
an eigenvariable
of an $(s2)$ with the index $\eta$.
\[
\infer[(s2)]{\Gam,\fal Y^{\eta}G}{\Gam,\fal Y^{\eta}G,G(V^{V^{\eta}})}
\]

\bdf\label{df:28}
\benu
\item\label{df:28.1}
{\rm The \textit{degree} $\dg(A)<\ome+\ome$.}
\[
\dg(A)=\left\{
\begin{array}{ll}
gr(A) & \mbox{{\rm if }} Gr(A)=0
\\
\ome+(Gr(A)-1) & \mbox{{\rm otherwise}}
\end{array}
\right.
\]

\item\label{df:28.2}
{\rm The \textit{height} $h(\Gam)=h(\Gam;P)$ of a sequent $\Gam$ in a preproof $P$.}
 \benu
 \item\label{df:28.21}
 $h(\Gam)=0$ {\rm if $\Gam$ is the end-sequent of $P$.}
 \item\label{df:28.22}
 $h(\Gam)=0$ {\rm if $\Gam$ is an upper sequent of a $(sub)$.}
 \item\label{df:28.23}
 $h(\Gam)=\ome$ {\rm if $\Gam$ is an upper sequent of a $(Q^{I}\mbox{{\rm -red}})$.
 
 In what follows assume that $\Gam$ is an upper sequent of a rule $J$ other than 
 $(sub), (Q^{I}\mbox{{\rm -red}})$ with the lower sequent $\Del$.}
 \item\label{df:28.24}
 $h(\Gam)=\max\{h(\Del),\dg(A)\}$ {\rm if $J$ is either a $(cut)$ with the cut formula $A$, or
 a $(BI)$ with the auxiliary formula $A$. }
 \item\label{df:28.25}
 $h(\Gam)=h(\Del)$ {\rm in other cases.}
 \eenu

\eenu
\edf

Relations between occurrences $A,B$ of formulas in a preproof 
such as `$A$ is a \textit{descendant} of $B$'
or equivalently `$B$ is an \textit{ancestor} of $A$',
and `an occurrence of inference rule is \textit{implicit} or \textit{explicit}' are defined as in \cite{Takeuti87,ptMahlo}.

\bdf\label{df:29}
{\rm Let $P$ be a preproof.}
\benu
\item
{\rm Let $\Del$ be a sequent in $P$. }
 \benu
 \item
 $\Del$ {\rm is in the \textit{explicit part} of $P$ if every rule below $\Del$ is either a explicit rule or
 a $(th)$, and $\Del$ is either an $(Ax)$ or a lower sequent of an explicit rule or a 
 $(th)$.}
 
 \item
 $\Del$ {\rm is a \textit{bar sequent} of $P$ if $\Del$ is not in the explicit part of $P$, and
 either $\Del$ is the end-sequent or an upper sequent of an explicit rule or a $(th)$ whose lower sequent is in the explicit part of $P$.}
 \eenu
 
 \item
 {\rm Let $\Del_{0}$ be a bar sequent of $P$. The \textit{end-piece} of $\Del_{0}$
 consists of the following sequents in $P$:
 $\Del_{0}$ is in the end-piece.
 If a lower sequent of a rule other than implicit rule is in the end-piece, then its upper sequents are in the end-piece.}
 
 \item
 {\rm An implicit rule is \textit{boundary rule} if its lower sequent is in an end-piece of $P$.}
 
 \item
 {\rm A triple $(J_{1},J_{2},J)$ of rules in $P$ is a \textit{suitable triangle}
 if $J_{i}$ is a boundary rule with its main formula $A_{i}$ for $i=1,2$,
 and $J$ is a $(cut)$:
 \[
 \infer[(cut)\, J]{\Gam,\Del}
 {
 \Gam,\lnot A
 &
 A,\Del
 }
 \]
where $\lnot A$ is a descendant of $A_{1}$, $A$ is a descendant of $A_{2}$ and
$\{\lnot A,A\}\cap(\Gam\cup\Del)=\emptyset$.

$A$ is said to be a \textit{suitable cut formula}.}
\eenu
\edf

\bprp\label{prp:10}
For a preproof $P$, $P$ contains no bar sequent iff
$P$ consists solely of explicit rules and $(th)$'s.
\eprp

In what follows a closed $s\in\cals$ is identified the ordinal $od(s)$.

\bdf\label{df:30}
{\rm Let $P$ be a preproof enjoying the pure variable condition.
A \textit{stack function} $sck$ for $P$ assigns an ordinal $sck(J)$ (the \textit{stack} of $J$)
to each occurrence $J$ of rules $(\exi^{I}\mbox{-red})$ and $(sub)$ in $P$.

Given a stack function $sck$, we assign 
ordinals $o(\Del)=o(\Del;P,sck), o(J)=o(J;P,sck)$ 
to each sequent $\Del$ and each line of a rule $J$ recursively as follows.
\benu
\item
$o(\Del)=1$ if $\Del$ is an $(Ax)$.

In what follows let $\Del$ be a lower sequent of a rule $J$ with upper sequents $\Gam$ (and $\Gam'$).

\item
$o(J)=o(\Gam)$ if $J$ is oner of rules $(Q\mbox{{\rm -red}}),(sub), (th)$.

\item
$o(J)=o(\Del)+ 1$ if $J$ is one of rules $(\lor), (\fal_{1}), (\exi_{1}),(c),(d), (s)$ and $(w)$.

\item
$o(J)=o(\Gam)\# o(\Gam')$ if $J$ is either a $(\land)$ or a $(cut)$.


  
\item
Let $J$ be a $(BI)$ with its main formula $\exi X F$.
\[
o(J)=\left\{
\begin{array}{ll}
\Ome_{s+1}\# o(\Gam) & \mbox{{\rm if }} \exi X F \mbox{ {\rm is stratified and }} s=st_{\Sig}(\exi X F)
\\
I\# o(\Gam) & \mbox{{\rm otherwise}}
\end{array}
\right.
\]

  






\item
If $J$ is a $(sub)^{s}$ of level $s$, then
$o(\Del)=\psi_{\Ome_{s+1}}(\gam\#\ome^{\alp})$ 
with $\gam=sck(J)$ and $\alp=o(J)=o(\Gam)$.

\item
If $J$ is an $(\exi^{I}\mbox{-red})$, then
$o(\Del)=\psi_{I}(\gam\#\ome^{\alp})$
with $\gam=sck(J)$ and $\alp=o(J)=o(\Gam)$.

\item
Let $J$ be a rule other than $(sub), (\exi^{I}\mbox{-red})$.
 \benu
 \item
 $o(\Del)=\ome_{m}(o(J))$ where $h(\Gam)=h(\Del)+m$ for an $m<\ome$.
  \item
 $o(\Del)=0$ if $h(\Del)<\ome\leq h(\Gam)$, cf.\,the condition in 
 Definition \ref{df:31}.\ref{df:31.1.5}.
 \eenu
 For ordinals $\alp$ and $m<\ome$, $\ome_{0}(\alp)=\alp$ and 
 $\ome_{m+1}(\alp)=\ome^{\ome_{m}(\alp)}$.
\eenu

Finally $o(P)=o(\Gam_{end};P,sck)$ for the end sequent $\Gam_{end}$ of $P$.

}
\edf

Note that we have $Gr(\exi X F)=0$ for the main formula $\exi X F$ of a $(BI)$, 
and $\exi X F$ is stratified if $\exi X F\not\in\Pi^{1}_{2}$,
cf\, Definition \ref{df:26}.\ref{df:26.BI}.
Then $od(st_{\Sig}(\exi X F))=I$ iff $\exi X F\in\Pi^{1}_{2}$ and a free variable $U^{s}$
occurs as a part of the formula $\exi X F$ such that $od(s)=I$,
while for an $s\in\cals$, $od(s)=I$ iff $s$ contains a free variable $V^{I}$ with the index $I$.
Suppose that the rule $(BI)$ is in a preproof with the pure variable condition and the
condition (\ref{eq:prp12}) in Proposition \ref{prp:12} is fulfilled for the preproof.
Then $od(st_{\Sig}(\exi X F))=I$ iff $\exi X F\in\Pi^{1}_{2}$ and a free variable $U^{s}$
occurs as a part of the formula $\exi X F$ such that
$s$ contains the eigenvariable $V^{I}$ of an $(s2)$ with type $I$.
\[
\infer[(s2)]{\Gam,\fal Y^{I} G}
{
 \infer*{\Gam,G(V^{V^{I}})}
 {
  \infer[(BI)]{\Del,\exi X F(U^{s(V^{I})})}{}
 }
}
 \]

\bprp\label{prp:11}
{\rm
Let $P$ be a preproof enjoying the pure variable condition.
Let $\Del$ be a sequent in $P$ and $Var(\Del)=\bigcup\{Var(A): A\in\Del\}$ for
the set $Var(A)$ of variables occurring in an index in the formula $A$.
Then the followings hold.
\benu

\item\label{prp:11b}
Each variable $U\in Var(\Del)$ is either an eigenvariable of a strong rule below $\Del$, or
$U\in Var(\Gam_{end})$ with the end-sequent $\Gam_{end}$ of $P$.

\item\label{prp:11c}
Let $U^{\eta}$ be a variable other than eigenvariables of strong rules in $P$,
and $s\in\cals$ such that $Var(s)=\emptyset$ and $s<I$.
For $P^{[s/U^{\eta}]}$, cf.\,Definition \ref{df:25},
 $P^{[s/U]}$ is a preproof enjoying the pure variable condition.
 
\item\label{prp:11d}
Let $U$ and $V$ be variables other than eigenvariables of strong rules in $P$.
Assume that $V$ does not occur as a part of a formula in $P$.
For each formula $A$ in $P$ let $A[V/U]$ denote the formula obtained from $A$ by replacing every occurrence of the variable $U$ as a part of a formula by the variable $V$.
$\Del[V/U]=\{A[V/U]: A\in\Del$, and $P[V/U]$ be the tree of sequents obtained from $P$ by replacing each sequent $\Del$ in $P$ by $\Del[V/U]$.

Then $P[V/U]$ is a preproof enjoying the pure variable condition.
\eenu
}
\eprp
\bprf

Proposition \ref{prp:11}.\ref{prp:11b} is seen inductively from below to above.

Propositions \ref{prp:11}.\ref{prp:11c} and \ref{prp:11}.\ref{prp:11d} are
 shown by induction on the depth of $P$ using Proposition \ref{prp:9}.

\eprf

\bprp\label{prp:12}
{\rm
Let $P$ be a preproof enjoying the pure variable condition.
Assume that $P$ satisfies the following condition:
\beqnarr
&&\mbox{The end-sequent of $P$ is a first-order sequent $\Gam_{end}$ such that}
\nonumber
\\
&&
\mbox{any $A\in\Gam_{end}$ is stratified and $st_{\Pi}(A)=0$}
\label{eq:prp12}
\eeqnarr
Let $\Del$ be a sequent in $P$.
\benu
\item\label{prp:12a}
If $h(\Del;P)<\ome$, then $\dg(A)<\ome$, i.e., $Gr(A)=0$ for any $A\in\Del$.

\item\label{prp:12b}
Let $U^{s}$ be a stratified variable occurring in $\Del$.
Then $I(s)\cap Fx<I$. In particular
 \benu
 \item\label{prp:12b1}
 if the main formula $A$ of a $(BI)$ is stratified, then $st_{\Sig}(A)<I$, and
 \item\label{prp:12b2}
  if the main formula $A$ of a $(w)$ is stratified, then $st_{\Pi}(A)<I$.
 \eenu
 
\item\label{prp:12c}
For any upper sequent $\Del$ of a $(sub)$, $Var(o(\Del))=\emptyset$ and $o(\Del)<I$.

\eenu
}
\eprp
\bprf
Proposition \ref{prp:12}.\ref{prp:12a} is seen inductively from below to above.
If $\Del$ is an upper sequent of s $(sub)^{s}$ of level $s<I$ and $A\in \Del$, then
$st_{\Pi}(A)\leq s$, and hence $Gr(A)=0$.

Proposition \ref{prp:12}.\ref{prp:12b} is shown inductively from below to above.
If $\Del$ is an upper sequent of a $(sub)^{s}$ of level $s<I$ with the eigenvariable $U$,
then $Var(s)=\emptyset$.
Hence $I(s)\cap F=I(s)<I$.

If $\Del$ is an upper sequent of a $(w)$ with the eigenvariable $U$ and 
the main formula $\fal X F$, then
$s=st_{\Pi}(\fal X F)$ and $Gr(\fal X F)=0$.
The assertion follows from IH.

If $\Del$ is an upper sequent of an $(s1)$ with the eigenvariable $U$ and 
the main formula $\fal X F$, then by IH,
we have $I(s)\cap F=I(\fal X F)\cap F<I$.

Proposition \ref{prp:12}.\ref{prp:12c} is seen from Proposition \ref{prp:11}.

\eprf

\bdf\label{df:31}
{\rm Let $P$ be a preproof enjoying the pure variable condition and $sck$ a stack function for $P$.
$P$ together with $sck$ 
is said to be a \textit{proof} (in ${\sf SBL}'$) if the following conditions are satisfied:
\benu

\item\label{df:31.2}
\beqnarr
\renewcommand{\theequation}{\ref{eq:prp12}}
&&\mbox{The end-sequent of $P$ is a first-order sequent $\Gam_{end}$ such that}
\nonumber
\\
&&
\mbox{any $A\in\Gam_{end}$ is stratified and $st_{\Pi}(A)=0$}
\eeqnarr
\addtocounter{equation}{-1}

\item\label{df:31.1.5}
Let $J$ be a rule with its lower sequent $\Del$ and an upper sequent $\Gam$ such that
$h(\Del)<\ome\leq h(\Gam)$.
Then the rule $J$ is a vacuous $(\exi^{I}\mbox{{\rm -red}})$.

Any rule $(Q^{I}\mbox{-red})\, J$ occurring in $P$ is in a series 
$(J_{0},\ldots,J_{n})$ of rules $(Q^{I}\mbox{-red})$, where 
$J=J_{i_{0}}$ for an $i_{0}\leq n$, each $J_{i+1}$ is immediately below $J_{i}$, 
there is a $k$ with $0\leq k\leq n$ such that each $J_{i}\,(i< k)$ is an $(\fal^{I}\mbox{-red})$,
while each $J_{i}\,(i\geq k)$ is an $(\exi^{I}\mbox{-red})$, and
there is no rule $(Q^{I}\mbox{-red})$ above $J_{0}$ nor below $J_{n}$.

\item\label{df:31.8}

Let $J$ be either an $(\exi^{I}\mbox{-red})$ or a $(sub)^{\mu}$,
$\Del$ the upper sequent of $J$, $\alp=o(\Del;P,sck)$ and
and $\gam=sck(J)$ the stack of $J$ with respect to the stack function $sck$.
Let $\sig=I$ when the rule is an $(\exi^{I}\mbox{{\rm -red}})$, and $\sig=\Ome_{\mu+1}$
when it is a rule $(sub)^{s}$.
Then for any index $s$ occurring above $J$

 \beqn\label{eq:df31.4a}
s\in\calh_{\gam}(\psi_{\sig}\gam)
\eeqn
and
\beqn\label{eq:df31.8}
\{\gam, \alp\}\subset\calh_{\gam}(\psi_{\sig}\gam)
\eeqn

where by an index $s$ occurring above $J$ we mean
\benu
\item\label{df:31.4a.1}
$T^{s_{0}}$ occurs above $J$ with $s\in I(s_{0})$, or
\item\label{df:31.4a.2}
$X^{s}$  occurs above $J$, or
\item\label{df:31.4a.3}
there is a rule $(\exi\mbox{-red})^{s}$ occurring above $J$, or
\item\label{df:31.4a.4}
there is a rule $(\exi\mbox{-red})\, J_{0}$ occurring above $J$ such that
$s=\psi_{I}(\alp\#\ome^{\bet})$ with $\alp=sck(J_{0})$ and $\bet=o(J_{0};P,sck)$.
\eenu

\item\label{df:31.4b}
Let $J$ be an $(\exi^{I}\mbox{-red})$ of type $\eta$ with the stack $\gam=sck(J)$,
$\Del$ the upper sequent of $J$ with $\alp=o(\Del;P,sck)$.
Then
\beqn\label{eq:df31.4b}
\eta\geq \psi_{I}(\gam\#\ome^{\alp})
\eeqn

\item\label{df:31.5}
every $(sub)$ is in an end-piece of a bar sequent.

\item\label{df:31.6}
the eigenvariable of a $(sub)$ does not occur in any explicit formula in the upper sequent of the 
$(sub)$.

\item\label{df:31.7}
each bar sequent $\Gam$ is the lower sequent of a vacuous $(sub)^{0}$ of level $0$.
The vacuous $(sub)$ is of the form
\[
\infer[(sub)^{0}]{\Gam}{\Gam}
\]
with an eigenvariable $U^{0}$ not occurring in $\Gam$.

\eenu
}
\edf


Clearly for any proof $P$, 
$o(P)<\Ome_{1}$.

For a first-order sequent $\Gam$ in the language $\calL$,
let $\Gam^{0}$ denote the sequent in $\calL'$ obtained from $\Gam$ by attaching the index $0$ to every second-order free variable and predicate constant occurring in $\Gam$.

\bprp\label{prp:13}
{\rm Let $\Gam$ be a first-order sequent $\Gam$ in $\calL$.
\benu
\item\label{prp:13a}
If $\Gam$ is derivable in ${\sf SBL}$, then so is $\Gam^{0}$ in ${\sf SBL}'$.

\item\label{prp:13b}
If there is a proof in ${\sf SBL}'$ ending with $\Gam^{0}$ and containing no bar sequent,
then $\Gam$ is (cut-free) derivable in the first-order sequent calculus ${\sf LK}$.
\eenu
}
\eprp
\bprf
Proposition \ref{prp:13}.\ref{prp:13a}.
Let $P$ be a ${\sf SBL}$-derivation of the first-order sequent $\Gam$.
We can assume that $P$ enjoys the pure variable condition,
$P$ contains no rule $(BI)_{2}$, cf. {\bf Remark} after Proposition \ref{prp:8}, and
any main formula $\exi X(A\subset X\subset B)$ of a $(\Pi^{1}_{2}\mbox{-Sep})$
is not $\Pi^{1}_{2}$, for
otherwise  $\exi X(A\subset X\subset B)\in\Pi^{1}_{2}\cap\Sig^{1}_{2}$, and
it is derivable by the rule $(BI)_{1}$ and $A\subset B$, i.e., from, e.g., $A\subset A\subset B$.

Then construct a proof $P^{0}$ of $\Gam^{0}$ from $P$ as follows:
attach the index $I$ to every undistinguished quantifiers,
attach the index $0$ to every predicate constant,
attach suitable indices to every second-order free variable from below to above.
Clearly the condition (\ref{eq:prp12}) in Definition \ref{df:31}.\ref{df:31.2} is enjoted,
and $st_{\Pi}(A)=0$ for any $A$ in the end-sequent $\Gam^{0}$, which is first-order.

In the resulting preproof $P_{0}$,
insert vacuous $(\exi^{I}\mbox{-red})$ immediately below a $(cut)$ such that
$h(\Del_{0},\Del_{1};P_{0})<\ome\leq h(\Del_{0},\lnot C;P_{0})$ for the lower sequent $\Del_{0},\Del_{1}$ and an upper sequent $\Del_{0},\lnot C$ of the $(cut)$.
Namely change 
\[
\infer[(cut)]{\Del_{0},\Del_{1}}
{
\Del_{0},\lnot C
&
C,\Del_{1}
}
\]
to
\[
\infer[(\exi^{I}\mbox{-red})]{\Del_{0},\Del_{1}}
{
 \infer[(cut)]{\Del_{0},\Del_{1}}
{
\Del_{0},\lnot C
&
C,\Del_{1}
}
}
\]
Then the condition in Definition \ref{df:31}.\ref{df:31.1.5} is fulfilled.
Note that $\dg(A)<\ome$, i.e., $Gr(A)=0$ for any $A\in\Del_{0}\cup\Del_{1}$ since the end-sequent 
$\Gam^{0}$
is first-order.
In particular no undistinguished quantifier $Q^{I}$ occurs in $\Del_{0}\cup\Del_{1}$.

Moreover insert vacuous $(sub)^{0}$ at bar sequents.
Note that any formula $B$ in any bar sequent is first-order, and hence $st_{\Pi}(B)=0$.

The resulting preproof is denoted $P^{0}$.
Any main formula of rules $(BI)$ and $(w)$ in $P^{0}$ is in $\Pi^{1}_{2}\cap\Sig^{1}_{2}$,
and $\eta=I$ for any main formula $\exi X^{\eta} F$ of rules $(c)$ in $P^{0}$.

A stack function $sck^{0}$ together with types of vacuous $(\exi^{I}\mbox{-red})$
is defined as follows.
First put $sck^{0}(J_{0})=0$ for any $(\exi^{I}\mbox{-red})\, J_{0}$.
Then the condition (\ref{eq:df31.4a}) is fulfilled
since any index $s$ occurring in $P^{0}$ is in $\{n,I+n:n<\ome\}$, and there is no rules 
$(\exi^{I}\mbox{-red})$ nor $(sub)$ in $P^{0}$.
Next the type $\eta$ of $(\exi^{I}\mbox{-red})\, J_{0}$ is defined to be
$\del_{0}=\psi_{I}(0\#\alp_{0})$ for $\alp_{0}=o(J_{0};P^{0},sck^{0})$.
Obviously $\alp_{0}\in\calh_{0}(0)\subset \calh_{0}(\psi_{I}0)$,
and the normality condition (\ref{eq:df31.8}) is fulfilled for $J_{0}$.
Then assign ordinals up to upper sequents of $(sub)^{0}\, J_{1}$.
Let $\alp_{1}=o(J_{1};P^{0},sck^{0})$, and pick an $n<\ome$ so that
$\alp_{0}<\ome_{n}(I+1)$ for any $(\exi^{I}\mbox{-red})\, J_{0}$ occurring above $J_{1}$
with $\alp_{0}=o(J_{0};P^{0},sck^{0})$.
Then let $sck^{0}(J_{1})=\ome_{n}(I+1)$, and
$\del_{1}=o(\Del;P^{0},sck^{0})=\psi_{\Ome_{1}}(\ome_{n}(I+1)\# \alp_{1})$
for the lower (bar) sequent $\Del$ of the $(sub)^{0}\, J_{1}$.
Then 
$\psi_{I}(0\#\alp_{0}),\ome_{n}(I+1)\#\alp_{1}\in\calh_{\ome_{n}(I+1)}(\psi_{\Ome_{1}}(\ome_{n}(I+1)))$.
Hence  the conditions (\ref{eq:df31.4a}) and (\ref{eq:df31.8}) are fulfilled for $J_{1}$.
 
Thus $P^{0}$ is a proof in ${\sf SBL}'$.

Proposition \ref{prp:13}.\ref{prp:13b} is seen from Proposition \ref{prp:10}.
Namely erase all the indices $0$ from the proof of $\Gam^{0}$ without bar sequent.
Then the result is a cut-free ${\sf LK}$-derivation of $\Gam$.
\eprf
\\

\noindent
By Propositions \ref{prp:8} and \ref{prp:13}, and the well-foundedness of $(\calh_{\veps_{I+1}}(0)\cap\Ome_{1},<)$
it suffices to show the following.

\bmlem\label{mlem}
For any proof $P$ and a stack function $sck$ in ${\sf SBL}'$ with a bar sequent,
we can construct another proof $P'$ and stack function $sck'$ 
with the same end-sequent such that
$o(P')<o(P)$.
\emlem

Main Lemma is proved in the next section \ref{sect:6}.

\section{Proof of Main Lemma \ref{mlem}}\label{sect:6}

Throughout this section $P$ together with a stack function $sck$ 
denotes a proof with a bar sequent.
For simplicity let us suppress stack functions in ordinals attached to sequents and rules.
Namely $o(\Gam;P,sck)$ [$o(J;P,sck)$] is denoted by
$o(\Gam;P)$ [$o(J;P)$], resp.

Each reduction, i.e., rewriting step is performed within the end-piece of a bar sequent $\Del_{0}$.
By Definition \ref{df:31}.\ref{df:31.7} the bar sequent $\Del_{0}$ is the lower sequent
of a vacuous $(sub)^{0}$ with its stack $\gam$,
\[
\infer[(sub)^{0}; \alp_{0}]{\Del_{b}; \psi_{\Ome_{1}}(\gam\#\ome^{\alp_{0}})}{\Del_{b};\alp_{0}}
\]
where and everywhere in this section,
$\Gam; \alp$ designates that $o(\Gam;P)=\alp$ for sequents $\Gam$ in $P$,
and $J;\alp$ that $o(J;P)=\alp$ for rules $J$ in $P$.
Also we see from 
(\ref{eq:prp12}) in Definition \ref{df:31}.\ref{df:31.2}
and the pure variable condition that 
each formula in an end-piece is stratified.

When $P$ is rewritten to another $P'$ below,
a stack function $sck'$ for $P'$ is defined in an obvious way
except otherwise stated explicitly.
Namely a rule $J'$ in $P'$ receives the same stack as one for the corresponding rule in $P$
in most cases.
In each step we need to verify that $P'$ is a proof and $o(P')<o(P)$.
In most cases this amounts to show that $P'$ together with a stack function $sck'$ fulfills the conditions
in (\ref{eq:df31.4a}), (\ref{eq:df31.8}) and (\ref{eq:df31.4b}).
\\

\noindent
{\bf Case 1}.
An explicit rule is in an end-piece of a bar sequent $\Del_{b}$ in $P$:

Let $J_{0}$ be one of the lowest explicit rule in the end-piece of $\Del_{b}$.
By (\ref{eq:prp12}) in Definition \ref{df:31}.\ref{df:31.2} 
the end-sequent of $P$ is a first-order sequent, and hence
$J_{0}$ is one of rules $(\land),(\lor),(\fal_{1}),(\exi_{1})$.
Consider the case when $J_{0}$ is a rule $(\fal_{1})$, and let $P$ be the following:
\[
\infer[J]{\Del_{b}; \psi_{\Ome_{1}}(\gam\#\ome^{\alp})}
{
 \infer*{\Del_{b};\alp}
{
 \infer[(\fal_{1})\, J_{0}]{\Gam;\bet+1}{\Gam,B(a);\bet}
}
}
\]
where $\fal x\, B(x)\in\Gam\cap\Del_{b}$ and $\Del_{b};\alp$ for the upper sequent
$\Del_{b}$ of the vacuous $(sub)^{0}\, J$ with its lower sequent $\Del_{b}$ and its stack $\gam=sck(J)$.
Note that by Definition \ref{df:31}.\ref{df:31.6} no $(sub)$ change explicit formulas,
and the end-piece ends with a vacuous $(sub)$ by Definition \ref{df:31}.\ref{df:31.7}.

Let $P'$ be the following.
\[
\infer[(\fal_{1})]{\Del_{b}; \psi_{\Ome_{1}}(\gam\#\ome^{\alp'})\# 1}
{
 \infer[J']{\Del_{b},B(a); \psi_{\Ome_{1}}(\gam\#\ome^{\alp'})}
 {
  \infer*{\Del_{b},B(a);\alp'}
  {
  \Gam,B(a);\bet
  }
  }
}
\]
We see from $\bet\ll\bet+1$ and Proposition \ref{prp:ll} that $\alp'\ll\alp$.
From this we see  $o(P')<o(P)$.
Let us verify that $P'$ is a proof.
The condition (\ref{eq:df31.4b}) on rules $(\exi^{I}\mbox{-red})$ in $P'$
is fulfilled by $\bet\ll\bet+1$.
We have $\alp\in\calh_{\gam}(\psi_{\Ome_{1}}\gam)$ by (\ref{eq:df31.8}) for $J$.
Hence $\alp'\in\calh_{\gam}(\psi_{\Ome_{1}}\gam)$ for the stack $\gam=sck'(J')$ of the vacuous 
$(sub)^{0}\, J'$ in $P'$.
Similarly we see that the conditions (\ref{eq:df31.4a}) and (\ref{eq:df31.8})  
on rules $(\exi^{I}\mbox{-red}), (sub)$ in $P'$ are fulfilled.
Therefore $P'$ is a proof.
\\

\noindent
{\bf Case 2}.
$\{\lnot A,A\}\subset\Del_{b}$ for a formula $A$ and a bar sequent $\Del_{b}$:

By (\ref{eq:prp12}) in Definition \ref{df:31}.\ref{df:31.2} $A$ is a first-order formula, and
$Gr(A)=0$.
\[
P=
\begin{array}{c}
\infer*{\Del_{b},\lnot A,A; \psi_{\Ome_{1}}(\gam\#\ome^{\alp})}{}
\end{array}
\msfiv
P'=
\begin{array}{c}
\infer[(Ax)]{\Del_{b},\lnot A,A;1}{}
\end{array}
\]
{\bf Case 3}.
The end-piece of a bar sequent $\Del_{b}$ contains a $(cut)$ of the following form:
\[
P=
\begin{array}{c}
\infer*{\Del_{b}}
{
 \infer[(cut);\alp\#\bet]{\Gam,\Del;\gam}
 {
  \Gam,\lnot A;\alp
  &
  A,\Del;\bet
  }
 }
\end{array}
\]
where $\lnot A\in\Gam\cup\Del$.
By Proposition \ref{prp:12}.\ref{prp:12a}
we have $h(\Gam,\Del)<\ome\Rarw h(\Gam,\lnot A)<\ome$.
In other words, $h(\Gam,\lnot A)=h(\Gam,\Del)+m$ for an $m<\ome$.
Thus $\gam=\ome_{m}(\alp\#\bet)$.
Let $P'$ be the following.
\[
P'=
\begin{array}{c}
\infer*{\Del_{b}}
{
 \infer[(th)]{\Gam,\Del;\alp'}
 {
  \Gam,\lnot A;\alp'
  }
 }
\end{array}
\]
Note that the height of an upper sequent of a $(sub)$ is defined to be $0$ in Definition \ref{df:28}.\ref{df:28.22},
and the height of an upper sequent of an $(\exi^{I}\mbox{-red})$ is equal to $\ome$ 
by Definition \ref{df:28}.\ref{df:28.23}.
Hence 
there is no $(sub)$ nor $(\exi^{I}\mbox{-red})$ in the height lowering part in $P'$.
Thus we see that $\alp'\ull\ome_{m}(\alp)\ll\gam$, and $P'$ is a proof such that $o(P')<o(P)$.
\\

\noindent
By virtue of {\bf Case 1-3}
we can assume that any end-piece of $P$ contains no explicit rule nor axiom.
Then we see as in Sublemma 12.9 of \cite{Takeuti87} that
$P$ contains a suitable triangle.

Before reducing suitable triangles, let us consider the following cases.
{\bf Cases 4-6} when a descendant of the main formula of a boundary rule $J_{0}=(d), (c)$
is changed by a $J_{1}=(\exi^{I}\mbox{-red})$,
and {\bf Cases 7-8} when a descendant of the main formula of a boundary rule $J_{0}=(s)$
is changed by a $J_{1}=(\fal^{I}\mbox{-red})$.
In each of these cases, $J_{0}$ and $J_{1}$ are exchanged.
When $J_{0}=(d)$, the distinguished rule $(d)$ is changed to a $(BI)$.
When $J_{0}=(s1)$, the strong rule $(s1)$ is changed to a weak rule $(w)$.
\\

\noindent
{\bf Case 4}.
A descendant of the main formula of a boundary rule $(d1)$ is changed by an 
$(\exi^{I}\mbox{-red})$:
Let $P$ be the following.
\[
P=
\begin{array}{c}
\infer[(\exi^{I}\mbox{-red})]{\Gam_{1},\exi X F'}
{
\infer*{\Gam_{0}, \exi X F' ; \del}
{
\infer[(\exi^{I}\mbox{-red})^{\eta}]{\Del',\exi X F'}
{
 \infer*{\Del,\exi X F;\del}
 {
  \infer[(d1)]{\Del_{0},\exi X F; \gam+1}
  {
   \Del_{0}, \exi X F, A\subset B;\gam
  }
 }
}
}
}
\end{array}
\]
where the lower rule $(\exi^{I}\mbox{-red})$ is a vacuous one such that 
$h(\Gam_{1},\exi X F')<\ome=h(\Gam_{0},\exi X F')$,
$F\equiv(A\subset X\subset B)$ with $(\exi X F)\in\Sig^{I}$ and $Gr(\exi X F)=1$.
Also
$F'\equiv F[\exi^{\eta}/\exi^{I}]\equiv(A'\subset X\subset B')\equiv(A[\fal^{\eta}/\fal^{I}]\subset X\subset B[\exi^{\eta}/\exi^{I}])$.
Hence $Gr(A')=Gr(B')=Gr(\exi X F')=0$, and $\exi X F'$ is stratified.
Note that there is no $(sub)^{s}$ 
between the boundary $(d1)$ and the $(\exi^{I}\mbox{-red})^{\eta}$
since the formula $\exi X F$ with $Gr(\exi X F)\neq 0$, and hence with
$st_{\Pi}(\exi X F)\geq I>s$  is not in the upper sequent of a $(sub)$.
All of these are seen from Definition \ref{df:26}.

Let $P'$ be the following.
\[
P'=
\begin{array}{c}
\infer[(BI)]{\Gam_{1},\exi X F'}
{
 \infer[(\land)]{\Gam_{1},\exi X F', F'(A')}
 {
  \infer*{\Gam_{1},\exi X F', A'\subset A';4}{}
  &
  \infer[(\exi^{I}\mbox{-red})]{\Gam_{1}, \exi X F', A'\subset B'}
  {
   \infer*{\Gam_{0},\exi X F', A'\subset B'; \del'}
   {
    \infer[(\exi^{I}\mbox{-red})^{\eta}]{\Del',\exi X F', A'\subset B'}
    {
     \infer*{\Del,\exi X F, A\subset B;\del'}
     {\Del_{0},\exi X F, A\subset B;\gam}
     }
    }
   }
  }
 }
\end{array}
\]
where $(A'\subset A')\equiv(\fal x[\lnot A'(x)\lor A'(x)])$ is derived from the axiom $\lnot A'(u),A'(u)$
with $Gr(A'(u))=0$ by two $(\lor)$'s followed by a $(\fal_{1})$.

It is easy to see that $o(\Gam,\exi X F,A\subset B;P')=\gam$ since
$Gr(F'(A'))=0$ and $dg(F'(A'))=gr(F'(A'))<\ome=h(\Gam_{0},\exi X F', A'\subset B';P')=h(\Gam_{0},\exi X F';P)$.
$\del'$ is an ordinal such that $\del'\ll \del$ by Proposition \ref{prp:ll}.
In particular $\psi_{I}(\alp\#\ome^{\del'})<\psi_{I}(\alp\#\ome^{\del})\leq\eta$ for the stack $\alp$ of 
the rules $(\exi^{I}\mbox{-red})^{\eta}$ by the condition (\ref{eq:df31.4b}).

Let $s$ be an index occurring in the formula $\exi X F$.
Then by the condition (\ref{eq:df31.4a}) we have
$s\in\calh_{\alp}(\psi_{I}\alp)$, and hence 
$s\in\calh_{\alp\#\ome^{\del}}(\psi_{I}(\alp\#\ome^{\del}))\cap I=\psi_{I}(\alp\#\ome^{\del})\leq\eta$.
Hence $st_{\Sig}(\exi X F')=\eta$.

Let us show $o(\Gam_{1},\exi X F';P')\ll o(\Gam_{1},\exi X F'; P)\,\{\eta\}$.
Let $\bet$ be the stack of 
the lower vacuous rule $(\exi^{I}\mbox{-red})$.
Then $o(\Gam_{1},\exi X F'; P)=\psi_{I}(\bet\#\ome^{\del})$, while
$o(\Gam_{1},\exi X F';P')=\ome_{m}(4\#\psi_{I}(\bet\#\ome^{\del'})\#\Ome_{\eta+1})$ 
for an $m<\ome$.

We see $o(\Gam_{1},\exi X F'; P')<o(\Gam_{1},\exi X F'; P)$ from
$\psi_{I}(\bet\#\ome^{\del'})\ll\psi_{I}(\bet\#\ome^{\del})$ and $\eta<\psi_{I}(\bet\#\ome^{\del})$.
The latter follows from (\ref{eq:df31.4a}), i.e.,
from $\eta\in\calh_{\bet}(\psi_{I}\bet)\cap I \subset\calh_{\bet\#\ome^{\del}}(\psi_{I}(\bet\#\ome^{\del}))\cap I=\psi_{I}(\bet\#\ome^{\del})$.
This yields $\Ome_{\eta+1}\ll\psi_{I}(\bet\#\ome^{\del})\,\{\eta\}$, and
$o(\Gam_{1},\exi X F';P')\ll o(\Gam_{1},\exi X F'; P)\,\{\eta\}$.
We see that the conditions (\ref{eq:df31.8}) and (\ref{eq:df31.4b}) 
is fulfilled for rules $(\exi^{I}\mbox{-red})$ in $P'$.

Consider the condition (\ref{eq:df31.4a}) for rules $(\exi^{I}\mbox{-red})$, e.g.,
for the lower vacuous rule $(\exi^{I}\mbox{-red})$.
There occur new indices, e.g., $\psi_{I}(\alp\#\ome^{\del'})$ for the rule $(\exi\mbox{-red})^{\eta}$
in $P'$, we need to show $\psi_{I}(\alp\#\ome^{\del'})\in\calh_{\bet}(\psi_{I}\bet)$.
We have $\psi_{I}(\alp\#\ome^{\del})\in\calh_{\bet}(\psi_{I}\bet)$ for the stack $\bet$ of the vacuous rule.
From $\psi_{I}(\alp\#\ome^{\del'})<\psi_{I}(\alp\#\ome^{\del})$, $\del'\ll\del$
and Proposition \ref{prp:ll}.\ref{prp:ll.6} we see that
$\psi_{I}(\alp\#\ome^{\del'})\ll\psi_{I}(\alp\#\ome^{\del})$.
In particular $\psi_{I}(\alp\#\ome^{\del'})\in\calh_{\bet}(\psi_{I}\bet)$.

Let $\Pi_{0}$ be an upper sequent of a $(sub)^{\mu}\, J$ occurring below $\Gam_{1},\exi X F'$ in $P$
with its lower sequent $\Pi_{1}$ and its stack $\bet=sck(J)$.
Assume $\alp_{0}' \ll \alp_{0}\, \{\eta\}$
for $\alp_{0}'=o(\Pi_{0};P')$ and $\alp_{0}=o(\Pi_{0};P)$,
Let the stack $\bet=stck'(J')$ of the corresponding rule $(sub)^{\mu}\, J'$ in $P'$.
We see that the condition (\ref{eq:df31.4a}) on $(sub)^{\mu}\, J'$ 
is fulfilled as above from  Proposition \ref{prp:ll}.\ref{prp:ll.6}.
For $\sig=\Ome_{\mu+1}$,
let $\alp_{1}=o(\Pi_{1};P)=\psi_{\sig}(\bet\#\ome^{\alp_{0}})$
 and $\alp_{1}'=o(\Pi_{1};P')=\psi_{\sig}(\bet\#\ome^{\alp_{0}'})$.
Then $\alp_{1}' \ll \alp_{1}\, \{\eta\}$ follows from Proposition \ref{prp:ll}.\ref{prp:ll.3} and (\ref{eq:df31.4a}),
$\eta\in\calh_{\bet}(\psi_{\sig}\bet)$.
Hence $o(P')\ll o(P)\,\{\eta\}$, and $o(P')<o(P)$.
\\

\noindent
{\bf Case 5}.
A descendant of the main formula of a boundary rule $(d2)$ is changed by an 
$(\exi^{I}\mbox{-red})$:
Let $P$ be the following.
\[
P=
\begin{array}{c}
\infer[(\exi^{I}\mbox{-red})]{\Gam_{1},\exi X F'}
{
 \infer*{\Gam_{0},\exi X F'}
 {
\infer[(\exi^{I}\mbox{-red})^{\eta}]{\Del',\exi X F'}
{
 \infer*{\Del,\exi X F;\del}
 {
  \infer[(d2)]{\Del_{0},\exi X F}
  {
   \Del_{0}, \exi X F, F(T^{s});\gam
  }
 }
}
}
}
\end{array}
\]
where the lower rule $(\exi^{I}\mbox{-red})$ is a vacuous one such that 
$h(\Gam_{1},\exi X F')<\ome=h(\Gam_{0},\exi X F')$,
$F'\equiv F[\exi^{\eta}/\exi^{I}]$, $Gr(\exi X F)\neq 0$, and $T^{s}$ is either
a predicate constant $R^{0}$ or a stratified free variable $U^{s}$ with $Var(s)=\emptyset$.
Also $o(\Del_{0},\exi X F;P)=\gam+1$ and $Gr(\exi X F')=0$ with stratified $\exi X F'$.
Similarly as in {\bf Case 4} we see that $st_{\Sig}(\exi X F')=\eta$,
and the following $P'$ is a proof such that $o(P')<o(P)$.
\[
P'=
\begin{array}{c}
\infer[(BI)]{\Gam_{1},\exi X F'}
{
 \infer[(\exi^{I}\mbox{-red})]{\Gam_{1},\exi X F',F'(T^{s})}
 {
  \infer*{\Gam_{0},\exi X F',F'(T^{s})}
  {
 \infer[(\exi^{I}\mbox{-red})^{\eta}]{\Del',\exi X F', F'(T^{s})}
 {
  \infer*{\Del,\exi X F, F(T^{s})}
  {
   \Del_{0},\exi X F, F(T^{s})
   }
  }
 }
 }
 }
\end{array}
\]
{\bf Case 6}.
A descendant of the main formula of a boundary rule $(c)$ is changed by an 
$(\exi^{I}\mbox{-red})$:
Let $P$ be the following.
\[
P=
\begin{array}{c}
\infer[(\exi^{I}\mbox{-red})]{\Gam_{1},\exi X^{\eta}F'}
{
 \infer*{\Gam_{0},\exi X^{\eta}F';\del}
 {
\infer[(\exi^{I}\mbox{-red})^{\eta}]{\Del',\exi X^{\eta}F'}
{
 \infer*{\Del,\exi X^{I}F}
 {
  \infer[(c)]{\Del_{0},\exi X^{I}F;\gam+1}
  {
   \Del_{0},\exi X^{I}F,F(T^{s});\gam
   }
  }
 }
 }
 }
\end{array}
\]
where the lower rule $(\exi^{I}\mbox{-red})$ is a vacuous one such that 
$h(\Gam_{1},\exi X F')<\ome=h(\Gam_{0},\exi X F')$,
$F'\equiv F[\exi^{\eta}/\exi^{I}]$.
There is no $(sub)$ between the boundary $(c)$ and the $(\exi^{I}\mbox{-red})$
since $Gr(\exi X^{I}F)\neq 0$ as in {\bf Case 4}.
By (\ref{eq:df31.4a}) and (\ref{eq:df31.4b})
we have $s\in\calh_{\gam}(\psi_{I}\gam)$ and $\psi_{I}\gam\leq\eta$ with the stack
of rules $(\exi^{I}\mbox{-red})^{\eta}$.
Hence $s<\eta$, and the rule $(c)$ in the following $P'$ is a legitimate one.


\[
P'=
\begin{array}{c}
\infer[(c)]{\Gam_{1},\exi X^{\eta}F'}
{
 \infer[(\exi^{I}\mbox{-red})]{\Gam_{1},\exi X^{\eta}F',F'(T^{s})}
 {
  \infer*{\Gam_{0},\exi X^{\eta}F',F'(T^{s});\del'}
  {
 \infer[(\exi^{I}\mbox{-red})^{\eta}]{\Del',\exi X^{\eta}F',F'(T^{s})}
 {
  \infer*{\Del,\exi X^{I}F,F(T^{s})}
  {
   \Del_{0},\exi X^{I}F,F(T^{s});\gam
   }
  }
 }
 }
 }
\end{array}
\]
We have $o(\Gam_{1},\exi X^{\eta}F';P)=\psi_{I}(\alp\#\ome^{\del})$
and $o(\Gam_{1},\exi X^{\eta}F';P')=\psi_{I}(\alp\#\ome^{\del'})+1$
for the stack $\alp$ of lower vacuous rules $(\exi^{I}\mbox{-red})$.
From $\del'\ll\del$ we see
$\psi_{I}(\alp\#\ome^{\del'})+1\ll\psi_{I}(\alp\#\ome^{\del})$.
We see easily that $P'$ is a proof such that $o(P')<o(P)$.
\\

\noindent
{\bf Case 7}. 
A descendant of the main formula of a boundary rule $(s1)$ is changed by an 
$(\fal^{I}\mbox{-red})$:
Let $P$ be the following.
\[
P=
\begin{array}{c}
\infer[(\exi^{I}\mbox{-red})]{\Gam_{1},\fal X F'}
{
 \infer*{\Gam_{0},\fal X F';\del}
 {
\infer[(\fal^{I}\mbox{-red})^{\eta}]{\Del',\fal X F'}
{
 \infer*{\Del,\fal X F;\del}
 {
  \infer[(s1)]{\Del_{0},\fal X F;\gam+1}
  {
   \infer*[P_{0}]{\Del_{0},\fal X F, F(U^{s});\gam}{}
   }
  }
 }
 }
 }
\end{array}
\]
where the lower rule $(\exi^{I}\mbox{-red})$ is a vacuous one such that 
$h(\Gam_{1},\fal X F')<\ome=h(\Gam_{0},\fal X F')$,
$F'\equiv F[\fal^{\eta}/\fal^{I}]$, $Gr(\fal X F)\neq 0$ and $s[I/U^{I}]=st_{\Pi}(\fal X F)$.
Since $\fal X F\in\Pi^{1}_{2}$ and $Gr(\fal X F)\neq 0$, $\fal^{I}$ occurs in $\fal X F$, i.e., 
$\fal X F\not\in\Sig^{I}$.
Therefore there occurs no $(sub)$ between the boundary $(s1)$ and the $(\fal^{I}\mbox{-red})$.
On the other, $Gr(\fal X F')=0$ and $\fal X F'$ is stratified.
Hence the rule $(w)$ in the following $P'$ is a legitimate one.
\[
P'=
\begin{array}{c}
\infer[(w)]{\Gam_{1},\fal X F'}
{
\infer[(\exi^{I}\mbox{-red})]{\Gam_{1},\fal X F',F'(U^{s'})}
{
 \infer*{\Gam_{0},\fal X F', F'(U^{s'});\del'}
 {
 \infer[(\fal^{I}\mbox{-red})^{\eta}]{\Del',\fal X F',F'(U^{s'})}
 { 
  \infer*{\Del,\fal X F,F(U^{s'});\del'}
  {
   \infer*[P_{0}']{\Del_{0},\fal X F, F(U^{s'});\gam'}{}
   }
  }
 }
 }
 }
\end{array}
\]
where $s'=st_{\Pi}(\fal X F')=s[\eta/U^{I}]$.
 $P_{0}'$ is obtained from $P_{0}$
by $P'_{0}=P_{0}^{[\eta/U^{I}]}$.

We have $\gam'\ull\gam\,\{\eta\}$, and $\gam'\ll\gam\,\{\eta\}$ if in $P_{0}$, there is a rule $(BI)$ 
with a main formula $\exi Y B$ such that the variable $U^{s}$ occurs as a part of $\exi Y B$, or
$s$ occurs in an index of a free variable in $\exi Y B$.
At such a rule $(BI)$, $I$ is added in $P$, 
while $\Ome_{\mu+1}$ is added in $P'$
for $\mu=st_{\Sig}((\exi Y B)[U^{s'}/U^{s}])$.

Let $\cals_{P_{0}}$ be the set of all indices $s_{0}$ such that 
either a free variable $V^{s_{0}}$ or a bound variable $Y^{s_{0}}$ occurs
in a main formula of a $(BI)$ in $P_{0}$.
Let $\eta'=\eta\#\bigcup\{I(s_{0}): s_{0}\in \cals_{P_{0}}\}$.
Then $\Ome_{\mu+1}\ll I\,\{\eta'\}$ for each such $(BI)$.
Hence $\del'\ll\del\,\{\eta'\}$.
On the other hand we have $\{\alp,\del,\eta'\}\subset\calh_{\alp}(\psi_{I}\alp)$ for
the stack $\alp$ of the lower vacuous rules $(\exi^{I}\mbox{-red})$ by 
(\ref{eq:df31.8}) and Definitions \ref{df:31}.\ref{df:31.4a.1}, \ref{df:31}.\ref{df:31.4a.2}.
Hence by Proposition \ref{prp:ll}.\ref{prp:ll.3}
we obtain $\psi_{I}(\alp\#\ome^{\del'})\ll\psi_{I}(\alp\#\ome^{\del})\,\{\eta'\}$ for
$o(\Gam_{1},\fal X F';P)=\psi_{I}(\alp\#\ome^{\del})$
and 
$o(\Gam_{1},\fal X F';P')=\psi_{I}(\alp\#\ome^{\del'})+1$.
Thus $o(\Gam_{1},\fal X F';P')\ll o(\Gam_{1},\fal X F';P)\,\{\eta'\}$, and we obtain
$o(P')\ll o(P)\,\{\eta'\}$, and $o(P')<o(P)$.

Let us verify that $P'$ is a proof.
Although $\eta$ is a new index in the upper part of $\Del,\fal X F, F(U^{s'})$,
there is no rule $(\exi^{I}\mbox{-red})$ nor $(sub)$ in the part 
since there is no $(sub)$ above the boundary $(s1)$ by Definition \ref{df:31}.\ref{df:31.5}.
Hence the conditions (\ref{eq:df31.4a}) and (\ref{eq:df31.8}) are enjoyed for the upper part.
The condition in (\ref{eq:df31.4b}) is fulfilled as we saw above.
We see that the conditions (\ref{eq:df31.4a}) and (\ref{eq:df31.8}) are fulfilled below 
$(\fal^{I}\mbox{-red})$ in $P'$ by Proposition \ref{prp:ll}.\ref{prp:ll.6} and (\ref{eq:df31.4a}) for $P$.
\\

\noindent
{\bf Case 8}.
A descendant of the main formula of a boundary rule $(s2)$ is changed by an 
$(\fal^{I}\mbox{-red})$:
Let $P$ be the following.
\[
P=
\begin{array}{c}
\infer[(\exi^{I}\mbox{-red})]{\Gam_{1},\fal X^{\eta} F'}
{
 \infer*{\Gam_{0},\fal X^{\eta} F'}
 {
\infer[(\fal^{I}\mbox{-red})^{\eta}]{\Del',\fal X^{\eta}F'}
{
 \infer*{\Del,\fal X^{I}F}
 {
  \infer[(s2)]{\Del_{0},\fal X^{I}F;\gam+1}
  {
   \Del_{0},\fal X^{I}F,F(U^{U^{I}}); \gam
   }
  }
 }
 }
 }
\end{array}
\]
where the lower rule $(\exi^{I}\mbox{-red})$ is a vacuous one such that 
$h(\Gam_{1},\fal X^{\eta} F')<\ome=h(\Gam_{0},\fal X^{\eta} F')$,
$F'\equiv F[\fal^{\eta}/\fal^{I}]$.

Let $P'$ be the following.
\[
P'=
\begin{array}{c}
\infer[(s2)]{\Gam_{1},\fal X^{\eta}F'}
{
\infer[(\exi^{I}\mbox{-red})]{\Gam_{1},\fal X^{\eta} F', F'(U^{U^{\eta}})}
{
 \infer*{\Gam_{0},\fal X^{\eta} F', F'(U^{U^{\eta}})}
 {
 \infer[(\fal^{I}\mbox{-red})^{\eta}]{\Del',\fal X^{\eta}F', F'(U^{U^{\eta}})}
 {
  \infer*{\Del,\fal X^{I}F, F(U^{U^{\eta}})}
  {
   \Del_{0},\fal X^{I}F,F(U^{U^{\eta}});\gam'
   }
  }
 }
 }
 }
\end{array}
\]
In $P'$, the index $U^{I}$ is replaced by $U^{\eta}$.
As in {\bf Case 7} we see that $P'$ is a proof such that
$o(P')<o(P)$.
\\

In the following cases let us reduce suitable triangles $(J_{1},J_{2},J)$,
where descendants of main formulas of $J_{1}$ and $J_{2}$
are not changed by any rules $(Q^{I}\mbox{-red})$ by virtue of {\bf Cases 4-8}.
\\
{\bf Case 9}.
$J_{1}$ is an $(s1)$ and $J_{2}$ is a $(d)$:
Let $P$ be the following.
\[
P=
\begin{array}{c}
\infer[J_{4}]{\Phi; \psi_{I}(\alp_{4}\#\ome^{\del})}
{
 \infer*{}
 {
  \infer[J_{3}]{}
  {
   \infer*{\Pi;\del}
   {
    \infer[J]{\Del,\Gam;\alp\#\bet}
    {
     \infer*{\Del,\fal X\lnot F;\alp}
     {
      \infer[(s1)\,J_{1}]{\Del_{0},\fal X\lnot F}
      {
       \Del_{0},\fal X\lnot F(U^{s})
       }
      }
     &
     \infer*{\exi X F,\Gam;\bet}
     {
      \infer[(d)\,J_{2}]{\exi X F,\Gam_{0}}
      {
       G,\exi X F,\Gam_{0}
       }
     }
    }
   }
  }
 }
}
\end{array}
\]
where $Gr(\exi X F)=1$, i.e., $\dg(\exi X F)=\ome=h(\Pi)$, 
$J_{2}$ is either a $(d1)$ with $G\equiv(A\subset B)$ and $F\equiv(A\subset X\subset B)$,
or a $(d2)$ with $G\equiv F(T^{s'})$.
$\Pi$ denotes the upper sequent of the uppermost $(\exi^{I}\mbox{-red})\, J_{3}$ 
below $J$.
$\Phi$ denotes the lower sequent of the lowest vacuous rule $(\exi^{I}\mbox{-red})\, J_{4}$.
In other words $\Phi$ is the uppermost sequent below the $(cut)\, J$ such that $h(\Phi)<\ome$.
Let $\alp_{n}=sck(J_{n})$ be the stack of the rule $J_{n}$ for $n=3,4$.

Note that no $(sub)$ occurs between $J$ and $\Phi$
since the height of the upper sequent of a $(sub)$ is defined to be $0$, 
cf.\,Definition \ref{df:28}.\ref{df:28.22}.
Furthermore there is no $(sub)$ between the $(s1)\, J_{1}$ and $(cut)\, J$,
and no $(sub)$ between the $(d)\, J_{2}$ and $J$ since $Gr(\exi X F)\neq 0$.

Let $P'$ be the following.
\[
P'=
\begin{array}{c}
\infer[(cut)]{\Phi}
{
\infer[J_{41}]{\Phi,\fal X\lnot F'; \psi_{I}(\alp_{4}'\#\ome^{\del_{1}})}
 {
  \infer*{}
  {
   \infer[J_{31}]{}
   {
   \infer[(\fal^{I}\mbox{-red})^{\eta}]{\Pi,\fal X\lnot F'}
   {
    \infer*{\Pi,\fal X\lnot F; \del_{1}}
    {
     \infer[(th)]{\Del,\Gam,\fal X\lnot F}
     { 
      \infer*{\Del,\fal X\lnot F;\alp}{}
     }
    }
   }
  }
 }
 }
&
\infer[J_{42}]{\exi X F',\Phi; \psi_{I}(\alp_{4}'\#\ome^{\del_{2}})}
 {
  \infer*{}
  {
   \infer[J_{32}]{}
   {
    \infer[(\exi^{I}\mbox{-red})^{\eta}]{\exi X F',\Pi}
    {
     \infer*{\exi X F,\Pi; \del_{2}}
     {
      \infer[(th)]{\exi X F,\Del,\Gam}
      {
       \infer*{\exi X F,\Gam;\bet}{}
      }
     }
    }
   }
  }
 }
}
\end{array}
\]
where $F'\equiv F[\exi^{\eta}/\exi^{I}]$.
Hence $\dg(\exi X F')=gr(\exi X F')<\ome=h(\Pi;P)$.
Then $o(\Del,\fal X\lnot F;P)=\alp=o(\Del,\fal X\lnot F;P')$ and
$o(\exi X F,\Gam;P)=\bet=o(\exi X F,\Gam;P')$.

From $\alp,\bet \ll\alp\#\bet$ we see that 
\beqn\label{eq:case9.1}
\del_{1}, \del_{2}\ll\del
\eeqn

The stack of the new rule $(\exi^{I}\mbox{-red})^{\eta}$ is defined to be $\alp_{3}=sck(J_{3})$, and
 the type $\eta$ of the new rules $(\fal^{I}\mbox{-red})^{\eta}$ and of
$(\exi^{I}\mbox{-red})^{\eta}$ is defined to be 
$\eta=\psi_{I}(\alp_{3}\# \ome^{\del_{2}})$ with $\del_{2}=o(\exi X F,\Pi;P')$. 
Let us verify the conditions (\ref{eq:df31.4a}), (\ref{eq:df31.8}) and (\ref{eq:df31.4b})
for the new rule $(\exi^{I}\mbox{-red})^{\eta}$.
(\ref{eq:df31.4b}) is obvious.
 (\ref{eq:df31.4a}) inherits from one for $J_{3}$ in $P$.
We have 
\beqn\label{eq:case9.2}
\alp_{3},\del\in\calh_{\alp_{3}}(\psi_{I}\alp_{3})
\eeqn
 by (\ref{eq:df31.8}) for $J_{3}$.
The condition (\ref{eq:df31.8}), $\alp_{3},\del_{2}\in\calh_{\alp_{3}}(\psi_{I}\alp_{3})$ 
follows from this and (\ref{eq:case9.1}).

Next let us increase stacks of the rules $(\exi^{I}\mbox{-red})\, J_{3i}$ by $\ome^{\del_{2}}+1$.
The stack of the rules $(\exi^{I}\mbox{-red})\, J_{31}$ and of $J_{32}$
is defined to be $\alp_{3}'=sck'(J_{31})=sck'(J_{32})=\alp_{3}\#\ome^{\del_{2}}\#1$.
We see that the conditions (\ref{eq:df31.4a}), (\ref{eq:df31.8}) and (\ref{eq:df31.4b})
 are fulfilled for $J_{3i}$ with $i=1,2$ as follows.
The new index $\eta\in\calh_{\alp_{3}'}(\psi_{I}\alp_{3}')$ for (\ref{eq:df31.4a}).
This is seen from 
(\ref{eq:case9.2}), (\ref{eq:case9.1}) and 
$\alp_{3}\# \ome^{\del_{2}}<\alp_{3}'$.
For $i=1,2$, we see 
$\alp_{3}',\del_{i}\in\calh_{\alp_{3}'}(\psi_{I}\alp_{3}')$ and
$\psi_{I}(\alp_{3}'\#\ome^{\del_{i}})\ll \psi_{I}(\alp_{3}\#\ome^{\del})$
from (\ref{eq:case9.1}) and (\ref{eq:case9.2}).
Thus the conditions (\ref{eq:df31.8}) and (\ref{eq:df31.4b}) are enjoyed for rules $J_{3i}$.

Let $K$ be a rule $(\exi^{I}\mbox{-red})$ occurring below $J_{3}$ in $P$,
and $\gam=sck(K)$ its stack.
Then the stack $\gam'$ of the corresponding rules $K'$ in $P'$
is defined to be $sck'(K')=sck(K)=\gam$, and let $\del'=o(K';P')\in\{\del_{1},\del_{2}\}$.
In particular the stack $sck'(J_{41})=sck'(J_{42})=sck(J_{4})=\alp_{4}$ of  the rules 
$(\exi^{I}\mbox{-red})\, J_{41}$ and of $J_{42}$.
We obtain $\gam,\ome^{\del'}\in\calh_{\gam}(\psi_{I}\gam)$ and
$\psi_{I}(\gam\#\ome^{\del'})\ll \psi_{I}(\gam\#\ome^{\del})$
from (\ref{eq:case9.1}) and $\{\gam,\del\}\subset\calh_{\gam}(\psi_{I}\gam)$.
Thus the conditions (\ref{eq:df31.8}) and (\ref{eq:df31.4b}) are enjoyed for rules $K'$.

Consider (\ref{eq:df31.4a}) for $K'$.
Let $\mu$ be the type of $J_{3}$.
Then we have $\mu\in\calh_{\gam}(\psi_{I}\gam)$ by Definition \ref{df:31}.\ref{df:31.4a.3}
for $K$.
On the other hand we have $\eta<\psi_{I}\alp_{3}'<\psi_{I}(\alp_{3}\#\ome^{\del})\leq\mu$
by (\ref{eq:df31.4b}) 
for $J_{3}$.
Hence $\eta<\mu\in\calh_{\gam}(\psi_{I}\gam)\cap I=\psi_{I}\gam$
and $\eta\in\calh_{\gam}(\psi_{I}\gam)$.
Moreover we have $\psi_{I}(\alp_{3}'\#\ome^{\del'})\in\calh_{\gam}(\psi_{I}\gam)$
by $\psi_{I}(\alp_{3}'\#\ome^{\del'})<\psi_{I}(\alp_{3}\#\ome^{\del})$, (\ref{eq:case9.1}),
Proposition \ref{prp:ll}.\ref{prp:ll.6} and $\psi_{I}(\alp_{3}\#\ome^{\del})\in\calh_{\gam}(\psi_{I}\gam)$.
Thus (\ref{eq:df31.4a}) is fulfilled for $K'$.

For $i=1,2$, we obtain
$\psi_{I}(\alp_{4}\#\ome^{\del_{i}})\ll \psi_{I}(\alp_{4}\#\ome^{\del})$.
Hence 
$o(\Phi;P')=\ome_{m}(\psi_{I}(\alp_{4}\#\ome^{\del_{1}})\#\psi_{I}(\alp_{4}\#\ome^{\del_{2}}))\ll
o(\Phi;P)$ for an $m<\ome$.

Finally let $S$ be a $(sub)^{\nu}$ occurring below $\Phi$ in $P$,
and $S'$ be the corresponding rule in $P'$.
Let $\gam=sck'(K')=sck(K)$, and $\sig=\Ome_{\nu+1}$.
We have $\psi_{I}(\alp_{3}\#\ome^{\del})\in\calh_{\gam}(\psi_{\sig}\gam)$ by 
Definition \ref{df:31}.\ref{df:31.4a.4}
for $K$.
Proposition \ref{prp:ll}.\ref{prp:ll.6} yields $\psi_{I}(\alp_{3}'\#\ome^{\del_{i}})\in\calh_{\gam}(\psi_{\sig}\gam)$
for $i=1,2$, and $\eta\in\calh_{\gam}(\psi_{\sig}\gam)$ by (\ref{eq:case9.1}) and
$\eta<\psi_{I}(\alp_{3}'\#\ome^{\del_{i}})<\psi_{I}(\alp_{3}\#\ome^{\del})$.

Thus the conditions (\ref{eq:df31.4a}) and (\ref{eq:df31.8}) for $K'$ are seen from
$o(K';P')\ll o(K;P)$ using Proposition \ref{prp:ll}.\ref{prp:ll.6} as above.
\\

\noindent
{\bf Case 10}.
$J_{1}$ is a $(w)$ and $J_{2}$ is a $(BI)$:
Let $P$ be the following.
\[
P=\begin{array}{c}
\infer[(sub)\, J_{3}]{\Phi}
{
 \infer*{\Pi; \del}
 {
  \infer[J]{\Del,\Gam}
  { 
   \infer*{\Del,\fal X\lnot F}
   {
    \infer[(w)\, J_{1}]{\Del_{0},\fal X\lnot F; \bet+ 1}
    {
     \Del_{0},\fal X\lnot F,\lnot F(U^{s});\bet
     }
    }
   &
    \infer*{\exi X F,\Gam}
    {
     \infer[(BI)\, J_{2}]{\exi X F,\Gam_{0};\alp_{1}}
     {
      F(A),\exi X F,\Gam_{0};\alp
      }
     }
    }
   }
  }
\end{array}
\]
where $Gr(\exi X F)=0$, $s=st_{\Pi}(\fal X\lnot F)=st_{\Sig}(\exi X F)$, and let $\sig=\Ome_{s+1}$.
Then $o(J_{2};P)=\sig\#\alp$ and $\alp_{1}=\ome_{m}(\sig\#\alp)$ for $m$ such that
$h(F(A),\exi X F,\Gam_{0})=\max\{h(\exi X F,\Del_{0}),\dg(F(A))\}=h(\exi X F,\Del_{0})+m$.
We see $s<I$ from Proposition \ref{prp:12}.\ref{prp:12b}.
Also $\Pi$ denotes the upper sequent of the uppermost $(sub)\, J_{3}$ of level$\leq s$ below $J$.

Note that no $(sub)$ changes the descendants of $\fal X\lnot F$ nor of $\exi X F$
by the condition in Definition \ref{df:26}.\ref{df:26.sub4.2}.

From Proposition \ref{prp:9}.\ref{prp:9b4} we see that $s=st_{\Pi}(\lnot F(U^{s}))$, and
$st_{\Pi}(\exi X F)=s+1$.
Hence from the Definition \ref{df:26}.\ref{df:26.sub4.1} of the rule $(sub)$, we see that
the level $\nu$ of any $(sub)^{\nu}$ occurring between $\exi X F,\Gam_{0}$ and $\Pi$ is
larger than $s$, $\nu>s$.
In particular no eigenvariable of a $(sub)$ occurring between $\Del,\Gam$ and $\Pi$
occurs in $\lnot F$.

Let $P'$ be the following.
\[
\infer[J_{3}']{\Phi}
{
 \infer*{\Pi; \del''}
 {
  \infer{\Del,\Gam,\Pi}
  {
   \infer*{\Del,\fal X\lnot F}{}
   &
   \infer*{\exi X F,\Gam,\Pi}
   {
    \infer[(cut)\, J_{2}']{\exi X F,\Gam_{0},\Pi;\alp_{1}'}
    {
     \infer[(sub)^{s}]{\Pi,\lnot F(A); \psi_{\sig}(\gam\#\ome^{\del'})}
     {
      \infer*{\Pi,\lnot F(U^{s});\del'}
      {
       \infer{\Del,\Gam,\lnot F(U^{s})}
       {
        \infer*{\Del,\lnot F(U^{s}),\fal X\lnot F}
        {
         \infer*{\Del_{0},\fal X\lnot F,\lnot F(U^{s});\bet}{}
          }
        &
        \infer*{\exi X F,\Gam}{}
        }
       }
      }
     &
     F(A),\exi X F,\Gam_{0};\alp
    }
   }
  }
 }
}
\]
where $h(\Pi,\lnot F(U^{s});P')=h(\Pi;P)=0$ for the upper sequent $\Pi,\lnot F(U^{s})$
of the new $(sub)^{s}$, $h(F(A),\exi X F,\Gam_{0};P')=h(F(A),\exi X F,\Gam_{0};P)$.
The rules occurring above $\Pi,\lnot F(U^{s})$ in $P'$ receives the same stack of the
corresponding rule in $P$.
From $\bet=o(\Del_{0},\fal X\lnot F,\lnot F(U^{s});P')=o(\Del_{0},\fal X\lnot F,\lnot F(U^{s});P)\ll
o(\Del_{0},\fal X\lnot F;P)=\bet+1$, we see for $\del'=o(\Pi,\lnot F(U^{s});P')$ that
\beqn\label{eq:case10.1}
\del'\ll\del
\eeqn

The stack $\gam$ of the new $(sub)^{s}$ is defined to be the stack $\gam=sck(J_{3})$ 
of the $(sub)\, J_{3}$  in $P$.
Then $\{\gam,\del'\}\subset\calh_{\gam}(\psi_{\sig}\gam)$,
$o(J_{2}';P')=\psi_{\sig}(\gam\#\ome^{\del'})\#\alp$, and for 
$\alp_{1}'=\ome_{m}(\psi_{\sig}(\gam\#\ome^{\del'})\#\alp)$
\beqn\label{eq:case10.2}
\alp_{1}'\ll\alp_{1}\,\{\psi_{\sig}(\gam\#\ome^{\del'})\}
\eeqn

Rules occurring between $J_{2}'$ and $J_{3}'$ in $P'$ receive the same stacks of the corresponding rules in $P$.
Then the condition (\ref{eq:df31.8}) is enjoyed for these $(sub)$'s by (\ref{eq:case10.2}).
Note that the level $\nu$ of any $(sub)^{\nu}$ between $J_{2}'$ and $J_{3}'$ 
is higher than $s$, $\nu>s$.
Then $\psi_{\sig}(\gam\#\ome^{\del'})<\sig\leq\Ome_{\nu}$ and
 for the stack $\lam$ of such a rule $(sub)^{\nu}$,
\beqn\label{eq:case10.2a}
\psi_{\sig}(\gam\#\ome^{\del'})\in\calh_{\lam}(\psi_{\Ome_{\nu+1}}\lam)
\eeqn
Since no eigenvariable of a $(sub)$ between $\Del,\Gam$ and $\Pi$ occurs in $\lnot F$,
such $(sub)$ does not change the descendants of $\lnot F(U^{s})$.

The stack of the $(sub)\, J_{3}'$ 
is increased by $\ome^{\del'}\#1$, i.e.,
$\gam'=sck'(J_{3}')=\gam\#\ome^{\del'}\#1$.
Let $\tau=\Ome_{\mu+1}$ with the level $\mu$ of $J_{3}$.
Then we have $\{\gam,\del\}\subset\calh_{\gam}(\psi_{\tau}\gam)$, and
$\del'\in \calh_{\gam}(\psi_{\tau}\gam)$ by (\ref{eq:case10.1}).
Hence $\gam'\in\calh_{\gam'}(\psi_{\tau}\gam)$.
Next from (\ref{eq:case10.2}), (\ref{eq:case10.2a})
 and Proposition \ref{prp:ll}.\ref{prp:ll.3}
we obtain 
\beqn\label{eq:case10.3}
\del''\ll\del\,\{\psi_{\sig}(\gam\#\ome^{\del'})\}
\eeqn
On the other hand we have $\sig\ll\alp_{1}$, and hence
$\sig\ll\del$ by Proposition \ref{prp:ll}.\ref{prp:ll.4}.
Hence
\beqn\label{eq:case10.3a}
\psi_{\sig}(\gam\#\ome^{\del'})\in\calh_{\gam'}(\psi_{\tau}\gam)
\eeqn
and
$\del''\in\calh_{\gam'}(\psi_{\tau}\gam)$ by (\ref{eq:case10.3}).
Therefore we obtain $\{\gam',\del''\}\subset\calh_{\gam'}(\psi_{\tau}\gam)$.
Thus the condition (\ref{eq:df31.8}) is enjoyed for
$J_{3}'$.

The condition (\ref{eq:df31.4a}) for $J_{3}'$ is enjoyed by (\ref{eq:case10.3a})
since no essentially new index occurs above $J_{3}'$.


Let us show 
$o(\Phi;P')=\psi_{\tau}(\gam'\#\ome^{\del''})\ll \psi_{\tau}(\gam\#\ome^{\del})=o(\Phi;P)$.
We have 
$\psi_{\tau}(\gam'\#\ome^{\del''})\ll \psi_{\tau}(\gam\#\ome^{\del})\, 
\{\psi_{\sig}(\gam\#\ome^{\del'})\}$
by (\ref{eq:case10.3}), Proposition \ref{prp:ll}.\ref{prp:ll.3} and (\ref{eq:case10.3a}).
Moreover we see from (\ref{eq:case10.1}) that
 for any $\alp>\gam\#\ome^{\del}$ and any $\rho$,
if $\{\gam,\del,\sig\}\subset\calh_{\alp}(\psi_{\rho}\alp)$,
then $\psi_{\sig}(\gam\#\ome^{\del'})\in\calh_{\alp}(\psi_{\rho}\alp)$.
$\psi_{\tau}(\gam'\#\ome^{\del''})\ll \psi_{\tau}(\gam\#\ome^{\del})$ is seen from 
Proposition \ref{prp:ll}.\ref{prp:ll.5}.

The stacks of $(sub)$'s below $J_{3}'$ remain the same.
(\ref{eq:df31.4a}) and (\ref{eq:df31.8}) are fulfilled for these $(sub)$'s by 
$\psi_{\tau}(\gam'\#\ome^{\del''})\ll \psi_{\tau}(\gam\#\ome^{\del})$.
\\

\noindent
{\bf Case 11}. 
$J_{1}$ is an $(s2)$ and $J_{2}$ is a $(c)$ with a main formula $\exi X^{I} F$:
Let $P$ be the following.
\[
P=\begin{array}{c}
\infer[J_{0};\alp]{\Phi;\bet}
{
 \infer*{}
 {
  \infer[J;\gam]{\Del,\Gam}
  {
   \infer*{\Del,\fal X^{I}\lnot F}
   {
    \infer[(s2)\,J_{1}]{\Del_{0},\fal X^{I}\lnot F;\del+1}
    {
     \infer*[P_{1}]{\Del_{0},\fal X^{I}\lnot F, \lnot F(U^{U^{I}});\del}{}
     }
    }
   &
   \infer*{\exi X^{I}F,\Gam}
   {
    \infer[(c)\,J_{2}]{\exi X^{I}F,\Gam_{0}; \xi+1}
    {
     F(V^{s}),\exi X^{I}F,\Gam_{0};\xi
     }
    }
  }
 }
}
\end{array}
\]
where $\bet=o(\Phi)$, $\alp=o(J_{0})$, $\gam=o(J)$, and 
$\del=o(\Del_{0},\fal X^{I}\lnot F,\lnot F(U^{U^{I}}))$,
and the lower sequent $\Phi$ of the rule $J_{0}$ 
denotes the uppermost sequent below $\Del,\fal X^{I}\lnot F$ such that
$h(\Phi)<h(\Del,\fal X^{I}\lnot F)$.

By Definition \ref{df:23}, $Gr(\fal X^{I}\lnot F)>1$, and hence
$h(\Del,\fal X^{I}\lnot F;P)>\ome$.
Since the height of upper sequents of $(Q^{I}\mbox{-red})$ is defined to be $\ome$,
we see that there is no $(Q^{I}\mbox{-red})$ between $J$ and $J_{0}$ in $P$.

Note that no $(sub)$ changes the descendants of $\fal X^{I}\lnot F$ nor of $\exi X^{I} F$
by the condition in Definition \ref{df:26}.\ref{df:26.sub4.2}.

From $\exi X^{I}F\not\in\Pi^{1}_{2}$, we see that
there is no $(sub)$ between $J_{2}$ and $J_{0}$ in $P$
since the height of the upper sequents of any $(sub)$ is defined to be $0$.


Let $P'$ be the following.
\[
\infer[(cut)]{\Phi;\bet'}
{
 \infer[J_{01};\alp_{1}]{\Phi,\lnot F(V^{s});\bet_{1}}
 {
  \infer*{}
  {
   \infer[J';\gam']{\Del,\Gam,\lnot F(V^{s})}
   {
    \infer*{\Del,\lnot F(V^{s}),\fal X^{I}\lnot F}
    {
     \infer*[P_{1}']{\Del_{0},\fal X^{I}\lnot F,\lnot F(V^{s});\del'}{}
     }
    &
    \infer*{\exi X^{I}F,\Gam}{}
    }
   }
  }
 &
 \infer[J_{02};\alp_{2}]{F(V^{s}),\Phi;\bet_{2}}
 {
  \infer*{}
  {
   \infer{F(V^{s}),\Del,\Gam}
   {
    \infer*{\Del,\fal X^{I}\lnot F}{}
    &
    \infer*{\exi X^{I}F,F(V^{s}),\Gam}
    {
     F(V^{s}),\exi X^{I} F,\Gam_{0};\xi
     }
    }
   }
  }
 }
\]
where $P_{1}'=(P_{1}^{[s/U^{I}]})[V/U]$, i.e., in $P_{1}$, 
replace first the occurrences of the variable $U^{I}$ in an index
by $s$, and then replace the occurrences of the variable $U$ 
as a part of formula by the variable $V$,
cf.\,Definition \ref{df:25}.

Let
$\bet'=o(\Phi;P')$, $\bet_{1}=o(\Phi,\lnot F(V^{s});P')$, $\bet_{2}=o(F(V^{s}),\Phi;P')$,
$\alp_{1}=o(J_{01};P')$, $\alp_{2}=o(J_{02};P')$.
Also $\del'=o(\Del_{0},\fal X^{I}\lnot F,\lnot F(V^{s});P')$ and $\gam'=o(J';P')$.

For $P'$ to be a proof, 
we need to verify the condition on rules $(sub)$ in Definition \ref{df:26}.\ref{df:26.sub4.1}, 
the condition (\ref{eq:df31.4b}) on rules $(\exi^{I}\mbox{-red})$, and
the conditions (\ref{eq:df31.4a}) and (\ref{eq:df31.8}) on rules $(sub),(\exi^{I}\mbox{-red})$.

First consider the condition on rules $(sub)$ in Definition \ref{df:26}.\ref{df:26.sub4.1}.
Since there is no $(sub)$ between $\exi X^{I} F,\Gam_{0}$ and $\Phi$ in $P$,
it suffices to examine a $(sub)$ occurring between $\Del_{0},\fal X^{I}\lnot F$ and $\Del,\fal X^{I}\lnot F$
with the added formula $\lnot F(V^{s})$ in $P'$.
From the same condition for the $(sub)$ in $P$ we see that $\fal X^{I}\lnot F\in\Pi^{1}_{2}$,
and hence $st_{\Pi}(\lnot F(V^{s}))\leq st_{\Pi}(\fal X^{I}\lnot F)$ by $s<I$.

Next consider the conditions on rules $(sub), (\exi^{I}\mbox{-red})$ in $P'$.
Let $K$ be a rule in $P$, which is either a $(sub)$ or an $(\exi^{I}\mbox{-red})$.
Assume that $K$ occurs either in $P_{1}$ or between $J_{1}$ and $J_{0}$.
We saw that $K$ is not between $J$ and $J_{0}$.
From $Gr(\fal X^{I}\lnot F)>1$ we see that $\fal X^{I}\lnot F$ is not in an upper sequent of a $(sub)$, 
which is in an end-piece.
Hence $K$ is not a $(sub)$.
Also from 
$h(\Del,\fal X^{I}\lnot F)>\ome$
and $h(\Del_{0},\fal X^{I}\lnot F,\lnot F(U^{U^{I}}))>\ome$,
we see that $K$ is not an $(\exi^{I}\mbox{-red})$.
Therefore there is no such rule $K$.

Let $\cals_{P_{1}}$ be the set of all indices $s_{1}$ such that 
either a free variable $W^{s_{1}}$ or a bound variable $Y^{s_{1}}$ occurs
in a main formula of a $(BI)$ in $P_{1}$.
Let $s'=s\#\bigcup\{I(s_{1}): s_{1}\in \cals_{P_{1}}\}$.

We have $\del'\ll\del\,\{s'\}$, and $\gam'\ll\gam\,\{s'\}$.
Hence $\alp_{1}\ll\alp\,\{s'\}$ and $\alp_{2}\ll\alp$.
Let $\Pi$ denote an upper sequent of $J_{0}$, and let $h=h(\Pi;P)$.
Then $h=h(\Phi)+m$ for an $m<\ome$.
From $\dg(F(V^{s}))<\dg(\exi X^{I}F)\leq h$ and $h(F(V^{s}),\Phi;P')<h$,
we see $\bet'\ll\bet\,\{s'\}$.

Let $K$ be the uppermost $(\exi\mbox{-red})$ below $J_{0}$ in $P$,
and $K'$ the corresponding rule in $K'$ with their stacks
$\gam=sck(K)=sck'(K')$.
Consider the conditions (\ref{eq:df31.8}) and (\ref{eq:df31.4b}) on $K'$.
We have $\alp_{K'}=o(K';P')\ll o(K;P)=\alp_{K}\,\{s'\}$, and
$s'\in\calh_{\gam}(\psi_{I}\gam)$ since the indices $s,s_{1}$ occur above $K$.
From $\{\alp_{K},\gam\}\subset\calh_{\gam}(\psi_{I}\gam)$ we obtain
$\alp_{K'}\in\calh_{\gam}(\psi_{I}\gam)$, and
$\psi_{I}(\gam\#\ome^{\alp_{K'}})<\psi_{I}(\gam\#\ome^{\alp_{K}})$.
Similarly we see that rules $(\exi^{I}\mbox{-red})$ below $J_{01}$ enjoy the conditions 
(\ref{eq:df31.4a}), (\ref{eq:df31.8}) and (\ref{eq:df31.4b}).

Next assume that $K$ is a $(sub)^{\mu}$ occurring below $J_{0}$, and consider the conditions 
(\ref{eq:df31.4a}) and (\ref{eq:df31.8}) on $K'$.
Let $\alp_{K}=o(K;P)$ and $\alp_{K'}=o(K';P')$.
We need to show that $\alp_{K'}\in\calh_{\gam}(\psi_{\sig}\gam)$,
where $\gam=sck(K)=sck'(K')$ and $\sig=\Ome_{\mu+1}$.
We have $\alp_{K'}\ll\alp_{K}\,\{s'\}$, and 
$s'\in\calh_{\gam}(\psi_{\sig}\gam)$ by (\ref{eq:df31.4a}) for $K$.
Then $\alp_{K'}\in\calh_{\gam}(\psi_{\sig}\gam)$.
\\

\noindent
{\bf Case 12}. 
$J_{1}$ is an $(s2)$ and $J_{2}$ is a $(c)$ with a main formula $\exi X^{\eta} F$ for an $\eta<I$:
Let $P$ be the following.
\[
P=\begin{array}{c}
\infer[J_{0};\alp]{\Phi;\bet}
{
 \infer*{}
 {
  \infer[J;\gam]{\Del,\Gam}
  {
   \infer*{\Del,\fal X^{\eta}\lnot F}
   {
    \infer[(s2)\,J_{1}]{\Del_{0},\fal X^{\eta}\lnot F;\del+1}
    {
     \infer*[P_{1}]{\Del_{0},\fal X^{\eta}\lnot F, \lnot F(U^{U^{\eta}});\del}{}
     }
    }
   &
   \infer*{\exi X^{\eta}F,\Gam}
   {
    \infer[(c)\,J_{2}]{\exi X^{\eta}F,\Gam_{0};\xi+1}
    {
     F(V^{s}),\exi X^{\eta}F,\Gam_{0};\xi
     }
    }
  }
 }
}
\end{array}
\]
where $\bet=o(\Phi)$, $\alp=o(J_{0})$, $\gam=o(J)$, and 
$\del=o(\Del_{0},\fal X^{\eta}\lnot F,\lnot F(U^{U^{\eta}}))$,
and the lower sequent $\Phi$ of the rule $J_{0}$ 
denotes the uppermost sequent below $\Del,\fal X^{\eta}\lnot F$ such that
$h(\Phi)<h(\Del,\fal X^{\eta}\lnot F)$.

Note that no $(sub)$ changes the descendants of $\fal X^{\eta}\lnot F$ nor of $\exi X^{\eta} F$
by the condition in Definition \ref{df:26}.\ref{df:26.sub4.2}.

From $\exi X^{\eta}F\not\in\Pi^{1}_{2}$, we see that
there is no $(sub)$ between $J_{2}$ and $J_{0}$ in $P$
since the height of the upper sequents of any $(sub)$ is defined to be $0$.

By Definition \ref{df:26}.\ref{df:26.c} and Proposition \ref{prp:12}.\ref{prp:12b}
we obtain $s<\eta$.
Then Proposition \ref{prp:9}.\ref{prp:9b3} with a limit $\eta$ yields
$st_{\Pi}(\lnot F(V^{s}))\leq st_{\Pi}(\fal X^{\eta}\lnot F)$ when $\fal X^{\eta}\lnot F\in\Pi^{1}_{2}$.

Let $P'$ be the following.
\[
\infer{\Phi;\bet'}
{
 \infer[J_{01};\alp_{1}]{\Phi,\lnot F(V^{s});\bet_{1}}
 {
  \infer*{}
  {
   \infer[J';\gam']{\Del,\Gam,\lnot F(V^{s})}
   {
    \infer*{\Del,\lnot F(V^{s}),\fal X^{\eta}\lnot F}
    {
     \infer*[P_{1}']{\Del_{0},\fal X^{\eta}\lnot F,\lnot F(V^{s});\del'}{}
     }
    &
    \infer*{\exi X^{\eta}F,\Gam}{}
    }
   }
  }
 &
 \infer[J_{02};\alp_{2}]{F(V^{s}),\Phi;\bet_{2}}
 {
  \infer*{}
  {
   \infer{F(V^{s}),\Del,\Gam}
   {
    \infer*{\Del,\fal X^{\eta}\lnot F}{}
    &
    \infer*{\exi X^{\eta}F,F(V^{s}),\Gam}
    {
     F(V^{s}),\exi X F,\Gam_{0};\xi
     }
    }
   }
  }
 }
\]
where $P_{1}'=(P_{1}^{[s/U^{\eta}]})[V/U]$, i.e., in $P_{1}$, 
replace first the occurrences of the variable $U^{\eta}$ in an index
by $s$, and then replace the occurrences of the variable $U$ 
as a part of formula by the variable $V$,
cf.\,Definition \ref{df:25}.

Let
$\bet'=o(\Phi;P')$, $\bet_{1}=o(\Phi,\lnot F(V^{s});P')$, $\bet_{2}=o(F(V^{s}),\Phi;P')$,
$\alp_{1}=o(J_{01};P')$, $\alp_{2}=o(J_{02};P')$.
Also $\del'=o(\Del_{0},\fal X^{\eta}\lnot F,\lnot F(V^{s});P')$ and $\gam'=o(J';P')$.

For $P'$ to be a proof, 
we need to verify the condition on rules $(sub)$ in Definition \ref{df:26}.\ref{df:26.sub4.1}, 
the condition (\ref{eq:df31.4b}) on rules $(\exi^{I}\mbox{-red})$, and
the conditions (\ref{eq:df31.4a}) and (\ref{eq:df31.8}) on rules $(sub),(\exi^{I}\mbox{-red})$.

We see that the condition on rules $(sub)$ in Definition \ref{df:26}.\ref{df:26.sub4.1} is fulfilled in $P'$
as in {\bf Case 11}.

Next consider the conditions on rules $(sub), (\exi^{I}\mbox{-red})$ in $P'$.
Let $K$ be a rule in $P$, which is either a $(sub)$ or an $(\exi^{I}\mbox{-red})$.
Assume that $K$ occurs either in $P_{1}$ or between $J_{1}$ and $J_{0}$.
Let $K'$ be the corresponding rule occurring in the left part of $J'$ in $P'$.
If the eigenvariable $U^{U^{\eta}}$ does not occur above $K$, then 
the new index $s$ does not occur above $K'$ except it occurs already above $K$,
and the ordinal remains the same.
In this case there is nothing to prove.
Assume that  $U^{U^{\eta}}$ occurs above $K$.

Let $\sig=I$ when $K$ is an $(\exi^{I}\mbox{-red})$, and $\sig=\Ome_{\mu+1}$
when $K$ is a $(sub)^{\mu}$.
Then $\eta<\sig$ is seen from $\eta<I$ when $K$ is an $(\exi^{I}\mbox{-red})$,
which is in $P_{1}$.
Also $\eta<\sig$ is seen from $\eta\leq st_{\Pi}(\fal X^{\eta}\lnot F)\leq \mu<\sig$
when $K$ is a $(sub)^{\mu}$, which is between $J_{1}$ and $J$, and
the formula $\fal X^{\eta}\lnot F$ is in the upper sequent of $K$, 
cf.\,Definition \ref{df:26}.\ref{df:26.sub4.1}.

Let $\gam=sck(K)=sck'(K')$ be the stack of the rule $K$ in $P$, 
and of the rule $K'$ in $P'$.
Let $\cals_{P_{1}}$ be the set of all indices $s_{1}$ such that 
either a free variable $W^{s_{1}}$ or a bound variable $Y^{s_{1}}$ occurs
in a main formula of a $(BI)$ in $P_{1}$.
Let $s'=s\#\bigcup\{I(s_{1}): s_{1}\in \cals_{P_{1}}\}$.
Since the variable $U^{U^{\eta}}$, i.e., the index $\eta$ as well as
indices $s_{1}$ in $P_{1}$ occurs above $K$,
we have 
$\{\eta\}\cup\bigcup\{I(s_{1}): s_{1}\in \cals_{P_{1}}\}\subset\calh_{\gam}(\psi_{\sig}\gam)$
by Definition \ref{df:31}.\ref{df:31.4a.1}.
Then $s<\eta\in\calh_{\gam}(\psi_{\sig}\gam)\cap\sig=\psi_{\sig}\gam$, and
\beqn\label{eq:case11.1}
s'\in\calh_{\gam}(\psi_{\sig}\gam)
\eeqn
Thus (\ref{eq:df31.4a}) is fulfilled for $K$.

Next let $\alp_{K}=o(K;P)$ and $\alp_{K'}=o(K';P')$.
We have $\alp_{K'}\ll\alp_{K}\,\{s'\}$.
By (\ref{eq:df31.8}) we have 
$\alp_{K}\in\calh_{\gam}(\psi_{\sig}\gam)$.
Hence $\alp_{K'}\in\calh_{\gam}(\psi_{\sig}\gam)$ by (\ref{eq:case11.1}).
Thus (\ref{eq:df31.8}) is fulfilled for $K'$.
(\ref{eq:df31.4b}) follows from  
$\alp_{K'}\in\calh_{\gam}(\psi_{I}\gam)$ and $\alp_{K'}<\alp_{K}$.


Finally let us show $o(P')<o(P)$.
We have $\del'\ull\del\, \{s'\}$.
Consider a $(sub)^{\mu}\, K$ occurring between $J_{1}$ and $J$.
Then $s<\eta\leq st_{\Pi}(\fal X^{\eta}\lnot F)\leq\mu<\Ome_{\mu+1}=\sig$,
and $\eta\in\calh_{\gam}(\psi_{\sig}\gam)\cap\sig$ with the stack $\gam=sck(K)$.
Hence $s'\in \calh_{\gam}(\psi_{\sig}\gam)$, and this yields 
$\gam'\ll\gam\, \{s'\}$.
We see that $\bet'\ll\bet\, \{s'\}$ as in {\bf Case 11}.
\\

\noindent
{\bf Case 13}.
The case when the suitable cut formula is a disjunction $A\lor B$.
\\
{\bf Case 14}.
The case when the suitable cut formula is an existential formula $\exi x A[x/u]$.

These cases are seen as in {\bf Case 11}.

This completes a proof of Main Lemma \ref{mlem}.


\end{document}